\documentclass{amsart}

\usepackage{amsmath}
\usepackage{amscd}

\usepackage{amsfonts,amssymb}

\theoremstyle{plain}
\newtheorem{thm}{Theorem}[section]
\newtheorem*{main}{Theorem}
\newtheorem{lemma}[thm]{Lemma}
\newtheorem{prop}[thm]{Proposition}
\newtheorem{cor}[thm]{Corollary}

\theoremstyle{remark}
\newtheorem{remark}[thm]{Remark}
\theoremstyle{definition}
\newtheorem{defn}[thm]{Definition}
\newtheorem{example}[thm]{Example}

\newcommand{\End}{\operatorname{End}}

\newcommand{\Vect}{{\sf{Vect}}}

\newcommand{\Aut}{\operatorname{Aut}}
\newcommand{\Emb}{\operatorname{Emb}}

\newcommand{\im}{\operatorname{im}}

\begin{document}
\title{Grassmannians of two-sided vector spaces}
\author{Adam Nyman}
\address{Department of Mathematics, University of Montana, Missoula, MT 59812-0864}
\email{NymanA@mso.umt.edu}
\keywords{Grassmannian, two-sided vector space, noncommutative vector bundle, bimodule}
\date{\today}
\thanks{2000 {\it Mathematics Subject Classification. } Primary 15A03, 14M15, 16D20; Secondary 14A22}
\thanks{The author was partially supported by the University of Montana University Grant program and by National Security Agency grant H98230-05-1-0021.}

\begin{abstract}
Let $k \subset K$ be an extension of fields, and let $A \subset M_{n}(K)$ be a $k$-algebra.  We study parameter spaces of $m$-dimensional subspaces of $K^{n}$ which are invariant under $A$.  The space $\mathbb{F}_{A}(m,n)$, whose $R$-rational points are $A$-invariant, free rank $m$ summands of $R^{n}$, is well known.  We construct a distinct parameter space, $\mathbb{G}_{A}(m,n)$, which is a fiber product of a Grassmannian and the projectivization of a vector space.  We then study the intersection $\mathbb{F}_{A}(m,n) \cap \mathbb{G}_{A}(m,n)$, which we denote by $\mathbb{H}_{A}(m,n)$.  Under suitable hypotheses on $A$, we construct affine open subschemes of $\mathbb{F}_{A}(m,n)$ and $\mathbb{H}_{A}(m,n)$ which cover their $K$-rational points.  We conclude by using $\mathbb{F}_{A}(m,n)$, $\mathbb{G}_{A}(m,n)$, and $\mathbb{H}_{A}(m,n)$ to construct parameter spaces of two-sided subspaces of two-sided vector spaces.
\end{abstract}

\maketitle

\tableofcontents

\section{Introduction}
Throughout this paper, $k \subset K$ is an extension of fields.  By a {\it two-sided vector space} we mean a $k$-central $K-K$-bimodule $V$ which is finite-dimensional as a left $K$-module.  Thus, a two-sided vector space on which $K$ acts centrally is just a finite-dimensional vector space over $K$.  The purpose of this paper is to continue the classification of two-sided vector spaces begun in \cite{pappacena} by constructing and studying parameter spaces of two-sided subspaces of $V$.

Instead of focusing exclusively on parameter spaces of two-sided subspaces of $V$, we take a more general perspective.  Let $A \subset M_{n}(K)$ be a $k$-algebra, and let $R$ be a $K$-algebra.  The functor $F_{A}(m,n):K-{\sf alg} \rightarrow {\sf Sets}$ defined on objects by
$$
F_{A}(m,n)(R)=\{ \mbox{free rank $m$ direct summands of $R^{n}$
which are $A$-invariant} \}
$$
and on morphisms by pullback is representable by a subscheme, $\mathbb{F}_{A}(m,n)$, of the Grassmannian of $m$-dimensional subspaces of $K^{n}$, $\mathbb{G}(m,n)$ \cite{rk}.  The scheme $\mathbb{F}_{A}(m,n)$ is related to two-sided vector spaces as follows.  Suppose $\phi:K \rightarrow M_{n}(K)$ is a $k$-central ring homomorphism and $K^{n}$ is made into a two-sided vector space, $K^{n}_{\phi}$, via $v \cdot x := v \phi(x)$.  Then the $K$-rational points of the scheme $\mathbb{F}_{\operatorname{im }\phi}(m,n)$ parameterize the two-sided $m$-dimensional subspaces of $K^{n}_{\phi}$.

There are other subschemes of $\mathbb{G}(m,n)$ which parameterize two-sided vector spaces as well.  In this paper, we study the geometry of $\mathbb{F}_{A}(m,n)$ and two other subschemes of $\mathbb{G}(m,n)$, $\mathbb{G}_{A}(m,n)$ and $\mathbb{H}_{A}(m,n)$, which have the same $K$-rational points as $\mathbb{F}_{A}(m,n)$.  Our justification for studying $\mathbb{G}_{A}(m,n)$ is that we are able to give a global description of it as an intersection of $\mathbb{G}(m,n)$ and the projectivization of a vector space (Theorem \ref{theorem.rep}).  Our justification for studying $\mathbb{H}_{A}(m,n)$ is that it is a subscheme of $\mathbb{G}_{A}(m,n)$ which has a smooth, reduced, irreducible open subscheme which covers its $K$-rational points (Theorem \ref{theorem.irr}).

We now describe $\mathbb{G}_{A}(m,n)$ by its functor of points, $G_{A}(m,n)$.  The $R$-rational points of this functor are the free rank $m$ summands of $R^{n}$, $M$, which have the property that the image of $\bigwedge^{m}M$ under the composition
$$
{\textstyle\bigwedge}^{m}R^{n} \overset{\cong}{\rightarrow} {\textstyle\bigwedge}^{m}(R \otimes_{K}K^{n}) \overset{\cong}{\rightarrow} R \otimes_{K}{\textstyle\bigwedge}^{m}K^{n}
$$
has an $R$-module generator of the form
$$
\underset{i}{\sum}r_{i} \otimes v_{i1} \wedge \cdots \wedge v_{im}
$$
where $\operatorname{Span }_{K}\{ v_{i1},\ldots,v_{im}\}$ is
$A$-invariant for all $i$ (see Definition \ref{def.phi}).  The
motivation behind this definition is that when $A-\{0\} \subset
GL_{n}(K)$ and $K^{n}$ is homogeneous as a $K \otimes_{k}A$-module
(see Section 2 for a description of the action of $K \otimes_{k}A$
on $K^{n}$), $G_{A}(m,n)$ solves the same parameterization problem
that $F_{A}(m,n)$ does, in the sense that
$G_{A}(m,n)(K)=F_{A}(m,n)(K)$.  Although it is not clear from the
definitions, the functors $G_{A}(m,n)$ and $F_{A}(m,n)$ are
distinct (Example \ref{example.notthesame}).

We prove in Section 3 that $G_{A}(m,n)$ has a simple global
description (Theorem \ref{theorem.rep}):

\begin{main}
Let
$$
{\textstyle\bigwedge}_{A}^{m}=\operatorname{Span}_{K}\{v_{1}\wedge \cdots \wedge v_{m}|\mbox{$v_{1},\ldots,v_{m}$ is a basis for an $A$-invariant subspace of $K^{n}$}\}.
$$
The functor $G_{A}(m,n)$ is represented by the pullback of the diagram
$$
\begin{CD}
& & \mathbb{P}_{K}(({\bigwedge}_{A}^{m})^{*}) \\
& & @VVV \\
\mathbb{G}(m,n) & \longrightarrow &
\mathbb{P}_{K}((\bigwedge^{m}K^{n})^{*})
\end{CD}
$$
whose horizontal is the canonical embedding, and whose vertical is induced by the inclusion ${\bigwedge}_{A}^{m} \rightarrow \bigwedge^{m}K^{n}$.
\end{main}
To the authors knowledge, there is no similar description of $\mathbb{F}_{A}(m,n)$.  The functorial description of ${\mathbb{G}}_{A}(m,n)$ allows us to describe the tangent space to ${\mathbb{G}}_{A}(m,n)$ (Theorem \ref{theorem.tangent}).

We define $\mathbb{H}_{A}(m,n)$ to be the pullback of the diagram
$$
\begin{CD}
& & \mathbb{G}_{A}(m,n) \\
& & @VVV \\
\mathbb{F}_{A}(m,n) & \longrightarrow & \mathbb{G}(m,n).
\end{CD}
$$
Suppose $S \subset K^{n}$ is a simple $K \otimes_{k}A$-module such that $\operatorname{dim }_{K}S=m$, and $K^{n}$ is $S$-homogeneous and semisimple as a $K \otimes_{k}A$-module.  In Section 4, we construct an affine open cover of the $K$-rational points of $\mathbb{F}_{A}(m,n)$.  Furthermore, when $K$ is infinite and $A$ is commutative, we construct an affine open cover of the $K$-rational points of $\mathbb{H}_{A}(m,n)$.  As a consequence, we prove the following (Theorem \ref{theorem.irr}):

\begin{main}
Suppose $K^{n} \cong S^{\oplus l}$ as $K \otimes_{k}A$-modules.  Then $\mathbb{F}_{A}(m,n)$ contains an open subscheme which is smooth, reduced, irreducible, of dimension $lm-m$ and has the same $K$-rational points as $\mathbb{F}_{A}(m,n)$.  Furthermore, if $K$ is infinite and $A$ is commutative, then $\mathbb{H}_{A}(m,n)$ contains an open subscheme which is smooth, reduced, irreducible, of dimension $lm-m$ and has the same $K$-rational points as $\mathbb{H}_{A}(m,n)$.
\end{main}

Now, let $V=K_{\phi}^{n}$ and let $W$ be a two-sided vector space.  In Section 6, we use $\mathbb{F}_{A}(m,n)$, $\mathbb{G}_{A}(m,n)$, and $\mathbb{H}_{A}(m,n)$ to construct three parameter spaces of two-sided subspaces of $V$ of rank $[W]$ (see Section 5 for the definition of rank).  We denote these parameter spaces by $\mathbb{F}_{\phi}([W],V)$, $\mathbb{G}_{\phi}([W],V)$, and $\mathbb{H}_{\phi}([W],V)$.  We then provide examples to show that, although $\mathbb{F}_{\phi}([W],V)$, $\mathbb{G}_{\phi}([W],V)$, and $\mathbb{H}_{\phi}([W],V)$ have the same $K$-rational points, $\mathbb{F}_{\phi}([W],V) \neq \mathbb{G}_{\phi}([W],V)$ and $\mathbb{F}_{\phi}([W],V) \neq \mathbb{H}_{\phi}([W],V)$ for certain $[W]$ and $V$.  As a consequence, $\mathbb{F}_{A}(m,n) \neq \mathbb{G}_{A}(m,n)$ and $\mathbb{F}_{A}(m,n) \neq \mathbb{H}_{A}(m,n)$ for certain $A$, $m$, and $n$.

We then show that, if $F$ is an extension field of $K$, then every element of $G_{\phi}([S],V)(F)$ and of $H_{\phi}([S],V)(F)$ is isomorphic to $F \otimes_{K} S$ as $F\otimes_{k}K$-modules (Theorem \ref{theorem.proj}).

We conclude by studying the geometry of the parameter spaces $\mathbb{F}_{\phi}([W],V)$, $\mathbb{G}_{\phi}([W],V)$, and $\mathbb{H}_{\phi}([W],V)$ in two cases.  In case $K/k$ is finite and Galois, we prove that the three spaces are equal to each other, and equal to the product of Grassmannians (Corollary \ref{cor.galois}).  In case $K$ is infinite, $\{ S_{i} \}_{i=1}^{r}$ consists of nonisomorphic simples with $\operatorname{dim }S_{i}=m_{i}$, and $V$ is semisimple with $l_{i}$ factors of $S_{i}$, we prove the following (Corollary \ref{cor.geometry1}):

\begin{main}
${\mathbb{F}}_{\phi}([S_{1}]+\cdots+[S_{r}],V)$ and ${\mathbb{H}}_{\phi}([S_{1}]+\cdots+[S_{r}],V)$ contain smooth, reduced, irreducible open subschemes of dimension $\sum_{i=1}^{r}l_{i}m_{i}-m_{i}$ which cover their $K$-rational points.
\end{main}

Aside from their significance as generalizations of Grassmannians, parameter spaces of two-sided subspaces of $V$, or Grassmannians of two-sided subspaces of $V$, are related to classification questions in noncommutative algebraic geometry.  The subject of noncommutative algebraic geometry is concerned, among other things, with classifying noncommutative projective surfaces (see \cite{stafford} for an introduction to this subject).  One important class of noncommutative surfaces, the class of noncommutative ruled surfaces, is constructed from noncommutative vector bundles.  Let $U$ and $X$ be schemes, and let $X$ be a $U$-scheme.  By a ``$U$-central noncommutative vector bundle over $X$", we mean an $\mathcal{O}_{U}$-central, coherent sheaf $X-X$-bimodule which is locally free on the right and left \cite[Definition 2.3, p. 440]{vdb}.

Two-sided vector spaces are related to noncommutative vector bundles as follows.  If $\mathcal{E}$ is a noncommutative vector bundle over an integral scheme $X$, $\mathcal{E}_{\eta}$ is a two-sided vector space over $k(X)$.  In addition, if $U=\operatorname{Spec }k$ and $X=\operatorname{Spec }K$, a $U$-central noncommutative vector bundle over $X$ is just a two-sided vector space.

Let $\mathcal{E}$ be a (commutative) vector bundle over $X$.  An important problem in algebraic geometry is to parameterize quotients of $\mathcal{E}$ with fixed Hilbert polynomial, and study the resulting parameter space.  We are interested in the analogous problem in noncommutative algebraic geometry: to parameterize $U$-central quotients of a $U$-central noncommutative vector bundle over $X$ with fixed invariants, and study the resulting parameter space.  Thus, the results in this paper address this problem when $U=\operatorname{Spec }k$ and $X=\operatorname{Spec }K$.

{\it Notation and conventions:}  We let ${\sf Sets}$ denote the category of sets and $K-{\sf alg}$ denote the category of commutative $K$-algebras.  For any scheme $Y$ over $\operatorname{Spec} K$, we let $h_{Y}$ denote the functor of points of $Y$, i.e. the functor $h_{Y}$ from the category $K-{\sf alg}$ to the category ${\sf Sets}$ is the functor $\operatorname{Hom}_{\operatorname{Spec} K}(\operatorname{Spec }-,Y)$.  Unless otherwise specified, all unlabeled isomorphisms are assumed to be canonical.  Finally, we suppose throughout that $A \subset M_{n}(K)$ is a $k$-algebra and $R$ is a commutative $K$-algebra.

Other notation and conventions will be introduced locally.

{\it Acknowledgments:}  I thank W. Adams, B. Huisgen-Zimmermann, C. Pappacena and N. Vonessen for helpful conversations, I thank A. Magidin for proving Lemma \ref{lemma.arturo}, and I thank S.P. Smith for a number of helpful comments regarding an earlier draft of this paper.

\section{Subfunctors of the Grassmannian}
Recall that the functor of points of the Grassmannian over $\operatorname{Spec }K$ is the functor $G(m,n):K-{\sf alg} \rightarrow {\sf Sets}$ defined on $R$ as the set of free rank $m$ summands of $R^{n}$, and defined on morphisms as the pullback \cite[Exercise VI-18, p. 261]{eisen}.

In this section, we define three subfunctors of $G(m,n)$, $F_{A}(m,n)$, $G_{A}(m,n)$, and $H_{A}(m,n)$.  We will see that $F_{A}(m,n)$ and $H_{A}(m,n)$ parameterize $m$-dimensional subspaces of $K^{n}$ which are invariant under the action of $A$, and $G_{A}(m,n)$ does so under suitable hypotheses on $A$.

Let $m=\sum_{i=1}^{n}r_{i}e_{i} \in R^{n}$, where $e_{i}$ is the standard unit vector.  We note that the action $r \otimes a \cdot m = \sum_{i=1}^{n}rr_{i}e_{i}a$ makes $R^{n}$ an $R \otimes_{k}A$-module.  We say that $M \subset R^{n}$ is {\it $A$-invariant} if $M$ is an $R \otimes_{k}A$-submodule.

\begin{defn} \label{definition.free}
Suppose $m$ is a nonnegative integer.  Let $F_{A}(m,n)(-): K-{\sf alg} \rightarrow {\sf Sets}$ denote the assignment defined on the object $R$ as
$$
\{\mbox{$M \in G(m,n)(R)|M$ is $A$-invariant}\}
$$
and on morphisms $\delta:R \rightarrow T$ as the pullback.  That is, $F_{A}(m,n)(\delta)(M)$ equals the image of the map
\begin{equation} \label{eqn.morph}
T \otimes_{R} M \rightarrow T \otimes_{R} R^{n} \overset{\cong}{\rightarrow} T^{n}
\end{equation}
whose left arrow is induced by inclusion $M \subset R^{n}$.
\end{defn}
The proof of the following result is straightforward, so we omit it.
\begin{lemma}
The assignment $F_{A}(m,n): K-{\sf alg} \rightarrow {\sf Sets}$ is a functor.
\end{lemma}
We call elements of $F_{A}(m,n)(R)$ {\it free rank $m$ $A$-invariant families over $\operatorname{Spec }R$}, or {\it free $A$-invariant families} when $m$, $n$ and $R$ are understood.

\begin{defn} \label{def.phi}
Let $M \subset R^{n}$ be a free rank $m$ summand.  We say $M$ is {\it generated by $A$-invariants over $\operatorname{Spec }R$} or is {\it generated by $A$-invariants} if $R$ is understood, if $\bigwedge^{m}M$ maps, under the composition
\begin{equation} \label{eqn.morph2}
\textstyle{\bigwedge^{m}M \rightarrow \bigwedge^{m}R^{n} \overset{\cong}{\rightarrow} \bigwedge^{m}(R \otimes_{K}K^{n}) \overset{\cong}{\rightarrow} R \otimes_{K} \bigwedge^{m}K^{n}}
\end{equation}
whose left arrow is induced by inclusion, to an $R$-module with generator of the form
\begin{equation} \label{eqn.wedge}
\underset{i}{\sum}r_{i} \otimes v_{i1}\wedge \cdots \wedge v_{im},
\end{equation}
where, for all $i$, $\{v_{i1},\ldots,v_{im}\}$ is a basis for a rank $m$ $A$-invariant subspace of $K^{n}$.
\end{defn}
For a discussion of the motivation behind this definition, see Remark \ref{remark.points}.

\begin{lemma} \label{lemma.functor}
Let $\delta:R \rightarrow T$ be a homomorphism of $K$-algebras, and let $M$ be a free rank $m$ summand of $R^{n}$ which is generated by $A$-invariants over $\operatorname{Spec }R$.  Then the image of $T \otimes_{R} M$ under (\ref{eqn.morph}) is a free rank $m$ summand of $T^{n}$ which is generated by $A$-invariants over $\operatorname{Spec }T$.
\end{lemma}

\begin{proof}
Suppose $M$ has basis $w_{1}, \ldots, w_{m} \in R^{n}$, and $w_{1} \wedge \cdots \wedge w_{m}$ maps to
$$
\underset{i}{\sum}r_{i} \otimes v_{i1}\wedge \cdots \wedge v_{im}
$$
under (\ref{eqn.morph2}).  Then $1 \otimes w_{1}, \ldots, 1 \otimes w_{m} \in T \otimes_{R}M$ are generators of $T \otimes_{R}M$, and, if we let $\overline{w_{1}}, \ldots, \overline{w_{m}}$ denote the images of $1 \otimes w_{1},\ldots, 1 \otimes w_{m}$ under (\ref{eqn.morph}), then $\overline{w_{1}},\ldots,\overline{w_{m}}$ generates the image of $T \otimes_{R} M$ under (\ref{eqn.morph}).  We claim $\overline{w_{1}} \wedge \cdots \wedge \overline{w_{m}} \in \bigwedge^{m}T^{n}$ maps to $\underset{i}{\sum}\delta(r_{i}) \otimes v_{i1}\wedge \cdots \wedge v_{im}$ under the composition
\begin{equation} \label{eqn.firstcomp}
\textstyle{\bigwedge^{m}T^{n} \overset{\cong}{\rightarrow} \bigwedge^{m}(T \otimes_{K}K^{n}) \overset{\cong}{\rightarrow} T \otimes_{K}\bigwedge^{m}K^{n}}.
\end{equation}
To prove the claim, we first note that a straightforward computation implies that
\begin{equation} \label{eqn.bigcommute}
\begin{CD}
\textstyle\bigwedge^{m}(T \otimes_{R}R^{n}) & \overset{\cong}{\rightarrow} & \textstyle\bigwedge^{m}T^{n} & \overset{\cong}{\rightarrow} & \textstyle\bigwedge^{m}(T \otimes_{K}K^{n}) & \overset{\cong}{\rightarrow} & T \otimes_{K}\textstyle\bigwedge^{m}K^{n} \\
@A{\cong}AA & & & & @AA{=}A \\
T \otimes \textstyle\bigwedge^{m}R^{n} & \overset{\cong}{\rightarrow} & T \otimes_{R} \textstyle\bigwedge^{m}(R \otimes_{R}K^{n}) & \overset{\cong}{\rightarrow} & T \otimes_{R} (R \otimes_{K} \textstyle\bigwedge^{m}K^{n}) & \overset{\cong}{\rightarrow} & T \otimes_{K}\textstyle\bigwedge^{m}K^{n}
\end{CD}
\end{equation}
commutes (recall our convention about unlabeled isomorphisms).  Furthermore, the image of $1 \otimes w_{1}\wedge \cdots \wedge w_{m} \in T \otimes_{R} \bigwedge^{m}R^{n}$ under the right-hand route of (\ref{eqn.bigcommute}) equals $\underset{i}{\sum}\delta(r_{i}) \otimes v_{i1}\wedge \cdots \wedge v_{im}$.  Therefore, the image of $1 \otimes w_{1}\wedge \cdots \wedge w_{m} \in T \otimes_{R} \bigwedge^{m}R^{n}$ under the left-hand route of (\ref{eqn.bigcommute}) equals $\underset{i}{\sum}\delta(r_{i}) \otimes v_{i1}\wedge \cdots \wedge v_{im}$.  Finally, the image of $1 \otimes w_{1}\wedge \cdots \wedge w_{m} \in T \otimes_{R} \bigwedge^{m}R^{n}$ under the first two maps of the left-hand route of (\ref{eqn.bigcommute}) equals $\overline{w_{1}} \wedge \cdots \wedge \overline{w_{m}} \in \bigwedge^{m}T^{n}$.  The claim, and hence the lemma, follows from the fact that the composition of the third and fourth arrows of the left-hand route of (\ref{eqn.bigcommute}) is the composition (\ref{eqn.firstcomp}).
\end{proof}

\begin{defn} \label{definition.grass}
Let $G_{A}(m,n)(-): K-{\sf alg} \rightarrow {\sf Sets}$ denote the assignment defined on the object $R$ as
$$
\{M \in G(m,n)(R)|M \mbox{ is generated by }A \mbox{-invariants}\}
$$
and on morphisms as the pullback.
\end{defn}
The next result follows immediately from Lemma \ref{lemma.functor}.

\begin{lemma}
$G_{A}(m,n):K-{\sf alg} \rightarrow {\sf Sets}$ is a functor.
\end{lemma}
We call elements of $G_{A}(m,n)(R)$ {\it free rank $m$ families generated by $A$-invariants over $\operatorname{Spec }R$}, or {\it free families generated by $A$-invariants} when $m$, $n$ and $R$ are understood.

\begin{remark} \label{remark.subset}
It follows immediately from Definition \ref{def.phi} that
$$
F_{A}(m,n)(K) \subset G_{A}(m,n)(K).
$$
\end{remark}
We now find sufficient conditions under which $F_{A}(m,n)(K) = G_{A}(m,n)(K)$.

\begin{lemma} \label{lemma.phiphiinvar}
Let $M$ be a free rank $m$ family generated by $A$-invariants over $\operatorname{Spec }R$.  If $M_{\mathfrak{m}}$ is $A$-invariant for every maximal ideal $\mathfrak{m}$ of $R$, then $M$ is $A$-invariant.  If $R$ is a field, $A-\{0\} \subset GL_{n}(K)$, and $K^{n}$ is homogeneous as a $K \otimes_{k} A$-module, then $M$ is $A$-invariant. \end{lemma}

\begin{proof}
Suppose $\mathfrak{m}$ is a maximal ideal of $R$, $a$ is an element of $A$, and $N = M a+M$.  The diagram
$$
\begin{CD}
R_{\mathfrak{m}} \otimes_{R} M & \rightarrow & R_{\mathfrak{m}} \otimes_{R} R^{n} & \overset{\cong}{\rightarrow} & R_{\mathfrak{m}}^{n} \\
@VVV @VV{=}V \\
R_{\mathfrak{m}} \otimes_{R} N & \rightarrow & R_{\mathfrak{m}} \otimes_{R} R^{n} & \overset{\cong}{\rightarrow} & R_{\mathfrak{m}}^{n}
\end{CD}
$$
whose left horizontals and left vertical are induced by inclusion, commutes, and the left horizontals are injective since localization is exact.  By Lemma \ref{lemma.functor}, the image, $\overline{M}$, of the top horizontal composition is generated by $A$-invariants.  Thus, by hypothesis, $\overline{M}$ is $A$-invariant.  Hence, the left vertical is surjective, and so the map
$$
M_{\mathfrak{m}} \rightarrow N_{\mathfrak{m}}
$$
induced by inclusion is an epimorphism.  It follows from \cite[Corollary 2.9, p. 68]{eisenbud} that $M=N$, and hence that $M$ is $A$-invariant.

Next, suppose $R$ is a field, $A-\{0\} \subset GL_{n}(K)$, $K^{n}$ is homogeneous as a $K \otimes_{k} A$-module, and $M$ has basis $w_{1},\ldots,w_{m} \in R^{n}$.  If $M$ were not $A$-invariant, then there would exist an $1 \leq i \leq m$ and an $a \in A$ such that $w_{i}a$ is not an element of $M$.  Thus, since $a$ is invertible, $w_{1}a \wedge \cdots \wedge w_{m} a$ would not be proportional to $w_{1}\wedge \cdots \wedge w_{m}$.

On the other hand, since $K^{n}$ is homogeneous, the determinant of the matrix corresponding to $a \in A$ acting on an $m$-dimensional, $A$-invariant, subspace $V$ is independent of $V$.  Thus, since $M$ is generated by $A$-invariants, $w_{1}a \wedge \cdots \wedge w_{m} a = c w_{1} \wedge \cdots \wedge w_{m}$ for some nonzero $c \in K$.  We conclude that $M$ is $A$-invariant.
\end{proof}

\begin{remark} \label{remark.points}
As a consequence of Lemma \ref{lemma.phiphiinvar}, if $A-\{0\} \subset GL_{n}(K)$, and $K^{n}$ is homogeneous as a $K \otimes_{k} A$-module,
$$
G_{A}(m,n)(K) = F_{A}(m,n)(K).
$$
Thus, these two functors parameterize the same object.  On the other hand, we will see in Theorem \ref{theorem.rep} that the scheme representing $G_{A}(m,n)$ has a simple global description as a pullback of $\mathbb{G}(m,n)$ and the projectivization of a vector space.  These two facts provide motivation for Definition \ref{def.phi}.
\end{remark}

\begin{defn} \label{definition.newgrass}
Let $H_{A}(m,n)(-): K-{\sf alg} \rightarrow {\sf Sets}$ denote the fibered product of functors $F_{A}(m,n) \times_{G(m,n)}G_{A}(m,n)$ induced by inclusion of $F_{A}(m,n)$ and $G_{A}(m,n)$ in $G(m,n)$ \cite[Definition VI-4, p. 254]{eisen}.
\end{defn}
We call elements of $H_{A}(m,n)(R)$ {\it free rank $m$ $A$-invariant families generated by $A$-invariants over $\operatorname{Spec }R$}, or {\it free $A$-invariant families generated by $A$-invariants} when $m$, $n$, and $R$ are understood.

\begin{remark} \label{remark.newpoints}
It follows from Remark \ref{remark.subset} that $H_{A}(m,n)(K)=F_{A}(m,n)(K)$.
\end{remark}

\section{Representability of $G_{A}(m,n)$ and $H_{A}(m,n)$}
It was proven in \cite{rk} that $F_{A}(m,n)$ is representable by a subscheme of the Grassmannian $\mathbb{G}(m,n)$.  The main result of this section is that $G_{A}(m,n)$ is representable by the intersection of $\mathbb{G}(m,n)$ and the projectivization of a vector space.  It will follow easily that $H_{A}(m,n)$ is representable as well.  We conclude the section by computing the tangent space to $G_{A}(m,n)$.

Let $\mathbb{P}(-)$ denote the projectivization functor.  That is, if $M$ is a $K$-module, we let $\mathbb{P}(M)$ denote the scheme whose $R$-rational points equal equivalence classes of epimorphisms $\tau: R \otimes_{K} M \rightarrow L$, where $L$ is an invertible $R$-module, such that $\tau_{1}:R \otimes_{K}M \rightarrow L_{1}$ is equivalent to $\tau_{2}:R \otimes_{K}M \rightarrow L_{2}$ iff there exists an isomorphism $\psi:L_{1} \rightarrow L_{2}$ such that $\tau_{2}=\psi \tau_{1}$.

Before proving that $G_{A}(m,n)$ is representable, we recall two preliminary facts.

\begin{lemma} \label{lemma.wedgedual}
Let $U$ be a subspace of $\bigwedge^{m}K^{n}$.  There is a natural
isomorphism
$$ 
R \otimes_{K} (U)^{*} \overset{\cong}{\longrightarrow}
(R\otimes_{K}U)^{*},
$$
and the canonical isomorphism ${\textstyle\bigwedge}^{m}R^{n}
\longrightarrow R\otimes_{K}{\textstyle\bigwedge}^{m}K^{n}$
induces an isomorphism
$$ 
(R\otimes_{K}{\textstyle\bigwedge}^{m}K^{n})^{*}
\overset{\cong}{\longrightarrow}
({\textstyle\bigwedge}^{m}R^{n})^{*}.
$$
\end{lemma}

We omit the straightforward proof of the next result.
\begin{lemma} \label{lemma.equal}
Let ${\sf F}$ denote the full subcategory of the category of $R$-modules consisting of finitely generated free $R$-modules.  Then the functor
$$
\operatorname{Hom}_{R}(-,R):{\sf F} \rightarrow {\sf F}
$$
is full and faithful.
\end{lemma}
We let
$$
{\textstyle\bigwedge}_{A}^{m}=\operatorname{Span}_{K}\{v_{1}\wedge \cdots \wedge v_{m}|\mbox{$v_{1},\ldots,v_{m}$ is a basis for an $A$-invariant subspace of $V$}\}.
$$
\begin{thm} \label{theorem.rep}
The functor $G_{A}(m,n)$ is represented by the pullback of the diagram
\begin{equation} \label{eqn.grass}
\begin{CD}
& & \mathbb{P}(({\bigwedge}_{A}^{m})^{*}) \\
& & @VVV \\
\mathbb{G}(m,n) & \longrightarrow & \mathbb{P}((\bigwedge^{m}K^{n})^{*})
\end{CD}
\end{equation}
whose horizontal is the canonical embedding, and whose vertical is induced by the inclusion ${\bigwedge}_{A}^{m} \rightarrow \bigwedge^{m}K^{n}$.
\end{thm}

\begin{proof}
By \cite[p. 260]{eisen}, it suffices to prove that the functor $G_{A}(m,n)$ is the pullback of functors
\begin{equation} \label{eqn.pullback0}
h_{\mathbb{G}(m,n)} \times_{h_{\mathbb{P}((\bigwedge^{m}K^{n})^{*})}} h_{\mathbb{P}(({\bigwedge}_{A}^{m})^{*})}
\end{equation}
induced by (\ref{eqn.grass}).

Let $M \subset R^{n}$ be a free rank $m$ summand.  We recall a preliminary fact.  The map $\bigwedge^{m}M \rightarrow \bigwedge^{m}R^{n}$ induced by the inclusion of $M$ in $R^{n}$ identifies $\bigwedge^{m}M$ with a summand of $\bigwedge^{m}R^{n}$.  Thus, the induced map
$$
\psi:({\textstyle\bigwedge}^{m}R^{n})^{*} \longrightarrow ({\textstyle\bigwedge}^{m}M)^{*}
$$
is an epimorphism.

We now prove that $G_{A}(m,n)$ equals (\ref{eqn.pullback0}).  We note that a free rank $m$ summand $M \subset R^{n}$ is an $R$-rational point of (\ref{eqn.pullback0}) iff the epimorphism
\begin{equation} \label{eqn.tau}
R \otimes_{K} ({\textstyle\bigwedge}^{m}K^{n})^{*} \overset{\cong}{\longrightarrow} ({\textstyle\bigwedge}^{m}R^{n})^{*} \overset{\psi}{\longrightarrow} ({\textstyle\bigwedge}^{m}M)^{*}
\end{equation}
whose left arrow is the map from Lemma \ref{lemma.wedgedual}, factors through the map
\begin{equation} \label{eqn.factor}
R \otimes_{K} ({\textstyle\bigwedge}^{m}K^{n})^{*} \longrightarrow R \otimes_{K} ({\textstyle\bigwedge}_{A}^{m})^{*}
\end{equation}
induced by the inclusion ${\bigwedge}_{A}^{m} \subset {\textstyle\bigwedge}^{m}K^{n}$.  Thus, to prove the result, it suffices to show that (\ref{eqn.tau}) factors through (\ref{eqn.factor}) iff $M$ is generated by $A$-invariants.  Now, (\ref{eqn.tau}) factors through (\ref{eqn.factor}) iff there exists a map $\gamma^{*}:(R \otimes_{K} {\bigwedge}_{A}^{m})^{*} \longrightarrow ({\bigwedge}^{m}M)^{*}$ making the diagram
\begin{equation} \label{eqn.gamma}
\begin{CD}
R \otimes_{K} ({\textstyle\bigwedge}^{m}K^{n})^{*} & & \overset{\cong}{\longrightarrow} & & (R \otimes_{K} {\textstyle\bigwedge}^{m}K^{n})^{*} & & \overset{\cong}{\longrightarrow}&  & ({\textstyle\bigwedge}^{m}R^{n})^{*} \\
@VVV & & @VVV & & @VV{\psi}V  \\
R \otimes_{K} ({\bigwedge}_{A}^{m})^{*} & & \overset{\cong}{\longrightarrow} & & (R \otimes_{K} {\bigwedge}_{A}^{m})^{*} & & \overset{\gamma^{*}}{\longrightarrow} & & ({\textstyle\bigwedge}^{m}M)^{*}
\end{CD}
\end{equation}
whose left and middle vertical are induced by inclusion, and whose top horizontals and bottom left horizontal are from Lemma \ref{lemma.wedgedual}, commute.  By Lemma \ref{lemma.wedgedual}, the left square of (\ref{eqn.gamma}) commutes.  Thus, by Lemma \ref{lemma.equal}, there exists a map $\gamma^{*}$ making (\ref{eqn.gamma}) commute iff there exists a map $\gamma:{\bigwedge}^{m}M \rightarrow R \otimes_{K} {\bigwedge}_{A}^{m}$ making the diagram
\begin{equation} \label{eqn.gamma2}
\begin{CD}
R \otimes_{K} {\textstyle\bigwedge}^{m} K^{n} & & \overset{\cong}{\longleftarrow} & & {\textstyle\bigwedge}^{m}R^{n} \\
@AAA & & @AAA \\
R \otimes_{K} {\bigwedge}_{A}^{m} & & \overset{\gamma}{\longleftarrow} & & {\textstyle\bigwedge}^{m}M
\end{CD}
\end{equation}
whose verticals are inclusions, commute.  This occurs iff $M$ is generated by $A$-invariants, i.e. iff $M \in G_{A}(m,n)(R)$.
\end{proof}
We denote the pullback of (\ref{eqn.grass}) by $\mathbb{G}_{A}(m,n)$.  The following result is now immediate:

\begin{cor}
$\mathbb{G}_{A}(m,n)$ is a projective subscheme of $\mathbb{G}(m,n)$.
\end{cor}
We also note that $H_{A}(m,n)$ is representable:

\begin{cor}
$H_{A}(m,n)$ is represented by the pullback of the diagram
$$
\begin{CD}
& & \mathbb{G}_{A}(m,n) \\
& & @VVV \\
\mathbb{F}_{A}(m,n) & \longrightarrow & \mathbb{G}(m,n).
\end{CD}
$$
whose arrows are induced by the inclusion of functors $G_{A}(m,n) \subset G(m,n)$ and $F_{A}(m,n) \subset G(m,n)$.
\end{cor}

\begin{proof}
Since $H_{A}(m,n)$ is defined as the fibered product $F_{A}(m,n) \times_{G(m,n)}G_{A}(m,n)$ induced by the inclusion of functors $G_{A}(m,n) \subset G(m,n)$ and $F_{A}(m,n) \subset G(m,n)$, the result follows immediately from \cite[p. 260]{eisen}.
\end{proof}

We end the section by computing the Zariski tangent space to $\mathbb{G}_{A}(m,n)$ at the $K$-rational point $E =\operatorname{Span }_{K}\{e_{1},\ldots, e_{m}\} \in G_{A}(m,n)(K)$.  Recall that if
$$
\Psi:G_{A}(m,n)(K[\epsilon]/(\epsilon^{2})) \rightarrow G_{A}(m,n)(K)
$$
is the map induced by the quotient $K[\epsilon]/(\epsilon^{2}) \rightarrow K$ sending $\epsilon$ to $0$, the Zariski tangent space to $\mathbb{G}_{A}(m,n)$ at the $K$-rational point $E \in G_{A}(m,n)(K)$ is the set
$$
T_{E}=\{M \in G_{A}(m,n)(K[\epsilon]/(\epsilon^{2}))|\Psi(M)=E\}
$$
with vector space structure defined as follows:  Suppose $\{f_{i}\}_{i=1}^{m}$ and $\{g_{i}\}_{i=1}^{m}$ are subsets of $K^{n}$.  If $M \in T_{E}$ has basis $\{e_{i}+\epsilon f_{i}\}_{i=1}^{m}$, $M' \in T_{E}$ has basis $\{e_{i}+\epsilon g_{i}\}_{i=1}^{m}$, and $a,b \in K$, we let $aM+bM' \in T_{E}$ denote the family with basis $\{e_{i}+\epsilon(af_{i}+bg_{i})\}$.  It is straightforward to check that the vector space structure is independent of choices made.

We define a map
$$
d:\operatorname{Hom}_{K}(E,K^{n}) \rightarrow \operatorname{Hom}_{K}(\textstyle\bigwedge^{m}E,\textstyle\bigwedge^{m}K^{n})
$$
as follows.  For $\psi \in \operatorname{Hom}_{K}(E,K^{n})$, we define $d(\psi)$ on totally decomposable wedges as
$$
d(\psi)(e_{1}\wedge \cdots \wedge e_{m})=\sum_{i=1}^{m}e_{1}\wedge \cdots \wedge e_{i-1} \wedge \psi(e_{i}) \wedge e_{i+1} \wedge \cdots \wedge e_{m}
$$
and extend linearly.  It is straightforward to check that $d$ is $K$-linear.

\begin{thm} \label{theorem.tangent}
Suppose $V=E \oplus L$ as a $K$-module for some $K$-submodule $L$ of $K^{n}$.  The tangent space to $\mathbb{G}_{A}(m,n)$ at $E \in \mathbb{G}_{A}(m,n)(K)$ is isomorphic to
$$
S_{E}=\{\psi \in \operatorname{Hom}_{K}(E,L) | \operatorname{im }d(\psi) \subset {\textstyle\bigwedge}_{A}^{m}\}.
$$
\end{thm}

\begin{proof}
We define a map
$$
\Phi:T_{E} \rightarrow S_{E}
$$
as follows:  let $M \in T_{E}$ have basis $\{e_{i}+\epsilon(s_{i} + t_{i})\}_{i=1}^{m}$ where $\{e_{1},\ldots,e_{m}\}$ is a basis for $E$, $s_{i} \in E$, and $t_{i} \in L$.  Define $\psi \in \operatorname{Hom}_{K}(E,K^{n})$ by $\psi(e_{i})=t_{i}$.  We let $\Phi(M)=\psi$, and we omit the straightforward proof of the fact that the definition of $\Phi$ is independent of choices made.

{\it Step 1:  We prove $\Phi$ is a well defined map of vector spaces.}  We omit the straightforward proof of the fact that as a map to $\operatorname{Hom}_{K}(E,K^{n})$, $\Phi$ is $K$-linear.  It remains to show that $\Phi(M) \in S_{E}$, i.e. that $\operatorname{im} d(\psi) \subset {\bigwedge}_{A}^{m}$.  Since $M$ is generated by $A$-invariants, $(e_{1}+\epsilon(s_{1}+t_{1}))\wedge \cdots \wedge (e_{m}+\epsilon(s_{m}+t_{m})) \in \bigwedge^{m}M$ maps, under (\ref{eqn.morph2}) to
\begin{equation} \label{eqn.tan}
\sum_{i}r_{i}\otimes v_{i_{1}}\wedge \cdots \wedge v_{i_{m}},
\end{equation}
where $r_{i} \in K[\epsilon]/(\epsilon)^{2}$ and $v_{i_{1}}\wedge \cdots \wedge v_{i_{m}} \in {\bigwedge}_{A}^{m}$.  Thus,
$$
(s_{1}+t_{1})\wedge e_{2} \wedge \cdots \wedge e_{m}+\cdots+e_{1}\wedge \cdots \wedge e_{m-1} \wedge (s_{m}+t_{m})= \sum_{j}a_{j}w_{j_{1}}\wedge \cdots \wedge w_{j_{m}},
$$
where $a_{j} \in K$ and $w_{j_{1}}\wedge \cdots \wedge w_{j_{m}} \in {\bigwedge}_{A}^{m}$.  This implies
\begin{equation} \label{eqn.tan2}
s_{1}\wedge e_{2} \wedge \cdots \wedge e_{m}+\cdots + e_{1}\wedge \cdots \wedge e_{m-1}\wedge s_{m} + t_{1}\wedge e_{2} \wedge \cdots \wedge e_{m}+\cdots + e_{1}\wedge \cdots \wedge e_{m-1}\wedge t_{m}
\end{equation}
is in ${\bigwedge}_{A}^{m}$.  Since $s_{i} \in E$, each of the first $m$ terms of (\ref{eqn.tan2}) is either $0$ or a multiple of $e_{1}\wedge \cdots \wedge e_{m}$, which is in ${\bigwedge}_{A}^{m}$.  Thus,
$$
t_{1}\wedge e_{2} \wedge \cdots \wedge e_{m}+\cdots + e_{1}\wedge \cdots \wedge e_{m-1}\wedge t_{m} \in {\textstyle\bigwedge}_{A}^{m},
$$
i.e. $d(\psi)(e_{1}\wedge \cdots \wedge e_{m}) \in {\bigwedge}_{A}^{m}$ as desired.

{\it Step 2:  We prove $\Phi$ is one-to-one and onto.}  If $\Phi(M)=0$ then $M$ has a basis $\{e_{i}\}_{i=1}^{m}$, and thus $M$ is the identity element of $T_{E}$.  This establishes the fact that $\Phi$ is one-to-one.

Let $\psi \in S_{E}$, and let $M \in T_{E}$ have basis $\{e_{i}+\epsilon\psi(e_{i})\}_{i=1}^{m}$.  To prove $\Phi$ is onto, we must prove that $M$ is generated by $A$-invariants.  By hypothesis,
$$
\psi(e_{1})\wedge e_{2} \wedge \cdots \wedge e_{m}+\cdots+e_{1}\wedge \cdots \wedge e_{m-1} \wedge \psi(e_{m})= \sum_{j}b_{j}u_{j_{1}}\wedge \cdots \wedge u_{j_{m}},
$$
where $b_{j} \in K$ and $u_{j_{1}}\wedge \cdots \wedge u_{j_{m}} \in {\bigwedge}_{A}^{m}$, and thus the image of $(e_{1}+\epsilon \psi(e_{1}))\wedge \cdots \wedge (e_{m}+\epsilon \psi(e_{m}))$ maps, under (\ref{eqn.morph2}), to
$$
\sum_{i}r_{i}v_{i_{1}}\wedge \cdots \wedge v_{i_{m}},
$$
where $r_{i} \in K[\epsilon]/(\epsilon)^{2}$ and $v_{i_{1}}\wedge \cdots \wedge v_{i_{m}} \in {\bigwedge}_{A}^{m}$.  Hence, $M \in T_{E}$, as desired.
\end{proof}

\section{Affine open subschemes of $\mathbb{F}_{A}(m,n)$ and $\mathbb{H}_{A}(m,n)$}
Suppose $S \subset K^{n}$ is a simple $K \otimes_{k}A$-module such that $\operatorname{dim }_{K}S=m$.  In this section, we assume $K^{n}$ is $S$-homogeneous and semisimple as a $K \otimes_{k}A$-module.  We study the subspace ${\bigwedge}_{A}^{m} \subset \bigwedge^{m}K^{n}$ in order to find conditions under which a free rank $m$ $A$-invariant family is generated by $A$-invariants.  We use our results in order to construct affine open subschemes of $\mathbb{F}_{A}(m,n)$ and $\mathbb{H}_{A}(m,n)$ which cover their $K$-rational points.

We suppose $K^{n}=S^{\oplus l}$, and we let $\pi_{i}:R^{lm} \rightarrow R^{m}$ denote projection onto the $(i-1)m+1$ through the $im$th coordinates.
\begin{lemma} \label{lemma.rankm}
Suppose $M \subset R^{lm}$ is $A$-invariant.  If $M$ is principally generated as an $R \otimes_{k}A$-module by $f$, and if $\pi_{i}(f) \in K^{m}$ for some $1 \leq i \leq l$, then $M$ is a free rank $m$ summand of $R^{lm}$, and $M$ is generated, as an $R$-module, by $f a_{1},\ldots,f a_{m}$ for some $a_{1},\ldots,a_{m} \in A$.
\end{lemma}

\begin{proof}
Suppose $\pi_{i}(f)=v \in K^{m}$.  Since $K^{n}$ is $S$-homogeneous, $S$ is simple, and $\operatorname{dim}_{K}S=m$, there exist $a_{1},\ldots,a_{m} \in A$ such that
$$
\{v a_{1},\ldots,v a_{m} \}
$$
is independent over $K$.  Thus, the $R$-module generated by $f a_{1},\ldots,f a_{m}$, which we denote by $\langle f a_{1},\ldots,f a_{m} \rangle$, is a free rank $m$ summand of $R^{lm}$.  To complete the proof of the lemma, it suffices to prove that $ \langle f a_{1},\ldots,f a_{m} \rangle$ is $A$-invariant.  To this end, we first prove $\pi_{i}|_{M}$ is injective.  Suppose $\pi_{i}(xf)=0$ for $x \in R \otimes_{k}A$.  Since $K^{n}=S^{\oplus l}$, $\pi_{i}(xf)=x\pi_{i}(f)$.  Thus, $xv=0$ in $R \otimes_{K}S$, so that $x \in \operatorname{ann }R \otimes_{K}S$.  Thus $x \in \operatorname{ann }R \otimes_{K}V$, again since $K^{n}=S^{\oplus l}$, so that $xf = 0$.

Now, suppose $a \in A$.  We prove $f a \in \langle f a_{1},\ldots,f a_{m} \rangle$.  For,
\begin{eqnarray*}
\pi_{i}(f a) & = & \pi_{i}(f) a \\
& = & v a \\
& = & b_{1}v a_{1}+\cdots+b_{m}v a_{m} \\
& = & \pi(b_{1}f a_{1}+\cdots b_{m}f a_{m}),
\end{eqnarray*}
where $b_{1},\ldots, b_{m} \in K$.  Since $\pi_{i}|_{M}$ is injective, we must have $f a=b_{1}f a_{1}+\cdots+b_{m}f a_{m}$, and the assertion follows.
\end{proof}
For the remainder of Section 4, we suppose $w_{i} \in S^{\oplus l}$ has one nonzero projection, $v_{i} \in S$, and $a_{1},\ldots,a_{m} \in A$ are such that
$$
\{v_{1} a_{1}, \ldots, v_{1} a_{m}\}
$$
is independent (such $a_{1},\ldots,a_{m}$ exist since $S$ is simple).

\begin{lemma} \label{lemma.prelim}
Suppose $A$ is commutative, $b_{1},\ldots,b_{r} \in K$, and consider the set
\begin{equation} \label{eqn.basis}
\{b_{1}w_{1} a_{1}+\cdots+b_{r}w_{r} a_{1},\ldots,b_{1}w_{1} a_{m}+\cdots+b_{r}w_{r} a_{m}\}.
\end{equation}
If (\ref{eqn.basis}) is nonzero, then (\ref{eqn.basis}) is a basis for an $A$-invariant subspace of $K^{n}$ of rank $m$.
\end{lemma}

\begin{proof}
We first note that, if $u \in S$ and $u \neq 0$, then $\{u a_{1},\ldots, u a_{m}\}$ is independent iff $\{v_{1} a_{1},\ldots,v_{1} a_{m}\}$ is independent.  For, since $S$ is a simple $K \otimes_{k}A$-module, there exists an $r \in K \otimes_{k}A$ such that $ru=v_{1}$.  Thus, since $K \otimes_{k}A$ is commutative, any dependence relation among $u a_{1},\ldots, u a_{m}$ is a dependence relation among $v_{1} a_{1},\ldots, v_{1} a_{m}$ and conversely.  Since we have assumed above that $\{v_{1} a_{1}, \ldots, v_{1} a_{m}\}$ is independent, we may conclude that if $u \in S$ is nonzero then $\{u a_{1},\ldots, u a_{m}\}$ is independent.

Suppose (\ref{eqn.basis}) is nonzero.  The fact that the set (\ref{eqn.basis}) is independent follows from the fact that some projection of (\ref{eqn.basis}) to a summand $S$ of $K^{n}=S^{\oplus l}$ is independent by the argument in the first paragraph.  To prove that the $K$-module generated by (\ref{eqn.basis}) is $A$-invariant, we note that the $K$-module generated by (\ref{eqn.basis}) is contained in the $K \otimes_{k}A$-module, $M$, generated by $b_{1}w_{1}  +\cdots+b_{r}w_{r}$.  On the other hand, by Lemma \ref{lemma.rankm}, $M$ is a free rank $m$ summand of $K^{ml}$.  Since the $K$-module generated by (\ref{eqn.basis}) is a free rank $m$ summand of $K^{n}$, it must equal $M$.
\end{proof}

\begin{prop} \label{proposition.phiinvar}
Let
$$
I=\{{\mathbf n}=(n_{1},\cdots,n_{r}) \in \mathbb{Z}_{\geq 0}^{r} | n_{1}+\cdots+n_{r}=m\},
$$
let $\{n_{1} \cdot 1,\ldots,n_{r}\cdot r\}$ denote the multiset with $n_{i}$ copies of $i$, and let
$$
w_{\mathbf{n}} = \underset{{\{(s_{1},\ldots,s_{m}) | \{s_{i}\}_{i=1}^{m}= \{n_{1}\cdot 1, \ldots, n_{r} \cdot r\}\}}}{\sum} w_{s_{1}} a_{1} \wedge \cdots \wedge w_{s_{m}} a_{m}.
$$
If $K$ is infinite, then $w_{{\bf n}}$ is an element of ${\bigwedge}_{A}^{m}$ for all ${\bf n} \in I$.
\end{prop}

\begin{proof}
If $r=1$, then $w_{\mathbf{n}}=w_{1}a_{1}\wedge \cdots \wedge w_{1}a_{m}$, and the result follows from Lemma \ref{lemma.prelim}.

Now suppose $r \geq 2$, so that $|I|=D=\binom{m+r-1}{m} \geq 2$.  We begin the proof of the case $r \geq 2$ with two preliminary observations.  First, each choice of $[(c_{1},\ldots,c_{r})]\in \mathbb{P}_{K}^{r-1}$ corresponds, via the $m$th Veronese map $\nu_{m}$, to a point in $\mathbb{P}_{K}^{D-1}$.  Since $K$ is infinite, no $D-2$ plane of $\mathbb{P}_{K}^{D-1}$ contains the image of $\nu_{m}$.

If ${\bf b} \in K^{r}$ with ${\bf b}=(b_{1},\ldots,b_{r})$, we let $b^{\bf n}=b_{1}^{n_{1}}\cdots b_{r}^{n_{r}}$.  Our second preliminary observation is that, for each $i_{0} \leq D$, there exist ${\bf b}_{1},\ldots,{\bf b}_{i_{0}} \in K^{r}$ with ${\bf b}_{i}=(b_{i1},\ldots,b_{ir})$, such that ${\{ (b_{i}^{\bf{n}})_{{\bf n}\in I}\}}_{i=1}^{i_{0}} \subset K^{D}$ is independent.  We prove this by induction on $i_{0} \geq 1$.  The case $i_{0}=1$ is trivial.  Now suppose the result holds for $1 \leq i_{0}$, where $i_{0} < D$.  Then there exist ${\bf b}_{1},\ldots,{\bf b}_{i_{0}} \subset K^{r}$ such that  ${\{ (b_{i}^{\bf{n}})_{{\bf n}\in I}\}}_{i=1}^{i_{0}} \subset K^{D}$ is independent.  Thus, since $\nu_{m}([{\bf b}_{i}])=[(b_{i}^{\bf{n}})_{{\bf n}\in I}]$, the subspace of $\mathbb{P}_{K}^{D-1}$ spanned by $\nu_{m}([{\bf b}_{1}]), \ldots, \nu_{m}([{\bf b}_{i_{0}}])$ is an $i_{0}-1$-plane.  Since $i_{0}<D$, the argument of the first paragraph implies there exists a ${\bf b}_{i_{0}+1} \subset K^{r}$ such that $\nu_{m}([{\bf b}_{i_{0}+1}])$ is not contained in the $i_{0}-1$-plane spanned by $\nu_{m}([{\bf b}_{1}]), \ldots, \nu_{m}([{\bf b}_{i_{0}}])$.  Thus, ${\{ (b_{i}^{\bf{n}})_{{\bf n}\in I}\}}_{i=1}^{i_{0}+1} \subset K^{D}$ is independent, as desired.  We conclude that there exist ${\bf b}_{1},\ldots,{\bf b}_{D} \in K^{r}$, such that ${\{ (b_{i}^{\bf{n}})_{{\bf n}\in I}\}}_{i=1}^{D} \subset K^{D}$ is independent.

We now prove the proposition.  For all $1 \leq i \leq D$, the vector
\begin{equation} \label{eqn.prelim0}
(b_{i1}w_{1} a_{1}+\cdots+b_{ir}w_{r} a_{1})\wedge \cdots \wedge (b_{i1}w_{1} a_{m}+\cdots+b_{ir}w_{r}a_{m})
\end{equation}
is $A$-invariant by Lemma \ref{lemma.prelim}.  Thus, by Remark \ref{remark.subset}, (\ref{eqn.prelim0}) is an element of ${\bigwedge}_{A}^{m}$.  In addition, (\ref{eqn.prelim0}) equals $\underset{{\bf n }\in I}{\sum}b_{i}^{{\mathbf{n}}}w_{{\mathbf{n}}}$.  Thus, it suffices to prove that $w_{\bf n} \in \operatorname{Span }\{\underset{{\bf n }\in I}{\sum}b_{i}^{{\mathbf{n}}}w_{{\mathbf{n}}}\}_{i=1}^{D}$ for all ${\bf n }\in I$.  To prove this, we note that since ${\{ (b_{i}^{\bf{n}})_{{\bf n}\in I}\}}_{i=1}^{D} \subset K^{D}$ is independent, $\{\underset{{\bf n }\in I}{\sum}b_{i}^{{\mathbf{n}}}w_{{\mathbf{n}}}\}_{i=1}^{D}$ is a set of $D$ independent vectors in $\operatorname{Span }\{w_{\bf n}| {\bf n} \in I\}$.  Since $|I|=D$, $\{\underset{{\bf n }\in I}{\sum}b_{i}^{{\mathbf{n}}}w_{{\mathbf{n}}}\}_{i=1}^{D}$ forms a basis for $\operatorname{Span }\{w_{\bf n}| {\bf n} \in I\}$.  Thus, $w_{\bf n} \in \operatorname{Span }\{\underset{{\bf n }\in I}{\sum}b_{i}^{{\mathbf{n}}}w_{{\mathbf{n}}}\}_{i=1}^{D}$ for all ${\bf n }\in I$, and the proof of the proposition follows.
\end{proof}

\begin{cor} \label{cor.before}
Suppose $K$ is infinite, $A$ is commutative, and $M \subset R^{lm}$ is $A$-invariant.  If $M$ is principally generated as an $R \otimes_{k}A$-module by $f$, and if $\pi_{i}(f) \in K^{m}$ for some $1 \leq i \leq l$, then $M$ is a free rank $m$ family generated by $A$-invariants over $\operatorname{Spec }R$.
\end{cor}

\begin{proof}
Throughout this proof, we let $[p]$ denote the set $\{1,\ldots,p\}$.  Let $\pi_{i}(f)=v \in K^{m}$.  By Lemma \ref{lemma.rankm}, there exist $a_{1},\ldots, a_{m} \in A$ such that $M$ is a free rank $m$ summand of $R^{lm}$ and $M$ is generated as an $R$-module by $f a_{1},\ldots,f a_{m}$.  Thus, $\bigwedge^{m}M$ is generated by $f a_{1}\wedge \cdots \wedge f a_{m}$ as an $R$-module.  On the other hand, $f=x_{1}u_{1}+\cdots+x_{l}u_{l}$ where $u_{i}=(0,\cdots,0,v_{i},0,\cdots,0) \in S^{\oplus l}$ has $i$th nonzero projection to $S$, and
$$
x_{i}=\sum_{j=1}^{n}r_{ij} \otimes b_{ij} \in R \otimes_{k}A.
$$
Thus, $\bigwedge^{m}M$ is generated by $(\sum_{i=1}^{l}x_{i}u_{i})a_{1}\wedge \cdots \wedge (\sum_{i=1}^{l}x_{i}u_{i})a_{m}$ which equals
$$
\underset{(i_{1},\ldots,i_{m}) \in [l]^{m}}{\sum} x_{i_{1}}u_{i_{1}} a_{1}\wedge \cdots \wedge x_{i_{m}}u_{i_{m}}a_{m}.
$$
Expanding further, we find the expression above equals
$$
\underset{(i_{1},\ldots,i_{m}) \in [l]^{m}}{\sum} \bigg{(}\underset{(j_{1},\ldots,j_{m}) \in [n]^{m}}{\sum} (r_{i_{1}j_{1}}\otimes b_{i_{1}j_{1}})\cdot u_{i_{1}} a_{1} \wedge \cdots \wedge (r_{i_{m}j_{m}}\otimes b_{i_{m}j_{m}}) \cdot u_{i_{m}} a_{m} \bigg{)},
$$
which equals
\begin{equation} \label{eqn.sum}
\underset{J}{\sum} r_{i_{1}j_{1}} \cdots r_{i_{m}j_{m}} u_{i_{1}} b_{i_{1}j_{1}} a_{1}\wedge \cdots \wedge  u_{i_{m}} b_{i_{m}j_{m}} a_{m}.
\end{equation}
where $J=([l]\times[n])^{m}$.  In order to prove $M$ is generated by $A$-invariants, it suffices to prove (\ref{eqn.sum}) is an element in the image of the composition
\begin{equation} \label{eqn.sum9}
R \otimes_{K} {\textstyle\bigwedge}_{A}^{m} \rightarrow R \otimes_{K} \textstyle{\bigwedge}^{m}K^{lm} \overset{\cong}{\rightarrow} \textstyle{\bigwedge}^{m}R^{lm},
\end{equation}
whose left arrow is induced by inclusion.

Let $S_{m}$ denote the $m$th symmetric group.  We note that $S_{m}$ acts on $J$ via $\sigma \cdot ((i_{1},j_{1}),\ldots,(i_{m},j_{m})) = ((i_{\sigma(1)},j_{\sigma(1)}),\ldots,(i_{\sigma(m)},j_{\sigma(m)}))$, and so $J$ is partitioned into the orbits of this action.  Thus, in order to prove (\ref{eqn.sum}) is an element in the image of (\ref{eqn.sum9}), it suffices to show
$$
w = \underset{\sigma \in \operatorname{S}_{m}}{\sum} u_{i_{\sigma(1)}}b_{i_{\sigma(1)}j_{\sigma(1)}} a_{1} \wedge \cdots \wedge u_{i_{\sigma(m)}}b_{i_{\sigma(m)}j_{\sigma(m)}} a_{m}
$$
is an element of ${\bigwedge}_{A}^{m}$, since (\ref{eqn.sum}) is an $R$-linear combination of images of terms of the form $1 \otimes_{K}w$ under (\ref{eqn.sum9}).  If we let $w_{q}=u_{i_{q}}b_{i_{q}j_{q}}$ for $1 \leq q \leq m$, and we let ${\bf n}=(1,\ldots,1) \in \mathbb{Z}_{\geq 0}^{m}$, $w$ is of the form $w_{\bf n}$ (defined in Proposition \ref{proposition.phiinvar}).  Since $K$ is infinite, the corollary follows from Proposition \ref{proposition.phiinvar}.
\end{proof}

We end this section by constructing collections of affine open subfunctors of $F_{A}(m,n)$ and $H_{A}(m,n)$ which cover their $K$-rational points.

For the remainder of this section, $B$ will denote the $K$-algebra $K[x_{1},\ldots,x_{lm-m}]$, and
$$
\langle (r_{1},\ldots,r_{lm}) \rangle \subset R^{lm}
$$
will denote the $R\otimes_{k}A$-submodule of $R^{lm}$ generated by $(r_{1},\ldots,r_{lm})$.  We will abuse notation as follows:  if $C$ and $D$ are $K$-algebras and $\psi \in h_{\operatorname{Spec }C}(D)$, we let $\psi:C \rightarrow D$ denote the induced map of rings.

For each $1 \leq i \leq l$, and each $R$, we define a map
$$
{\Phi_{i}}_{R}:h_{\operatorname{Spec }B}(R) \rightarrow G(m,n)(R)
$$
as follows:  if $\psi \in h_{\operatorname{Spec }B}(R)$, let
\begin{equation} \label{eqn.bigphi}
{\Phi_{i}}_{R}(\psi)=\langle (\psi(x_{1}),\ldots,\psi(x_{(i-1)m}),1,0,\ldots,0,\psi(x_{(i-1)m+1}),\ldots,\psi(x_{lm-m})) \rangle.
\end{equation}

\begin{lemma} \label{lemma.welldef}
${\Phi_{i}}_{R}$ is a well defined map of sets, and induces a natural transformation $\Phi_{i}:h_{\operatorname{Spec }B} \rightarrow G(m,n)$ which factors through the inclusion $F_{A}(m,n) \rightarrow G(m,n)$.  Furthermore, if $K$ is infinite and $A$ is commutative, then $\Phi_{i}$ factors through the inclusion $H_{A}(m,n) \rightarrow G(m,n)$.
\end{lemma}

\begin{proof}
Suppose $\psi \in h_{\operatorname{Spec }B}(R)$.  By Lemma \ref{lemma.rankm}, ${\Phi_{i}}_{R}(\psi)$ is a free rank $m$ $A$-invariant submodule of $R^{lm}$, whence the first assertion.  The proof that ${\Phi_{i}}_{R}$ induces a natural transformation $\Phi_{i}:h_{\operatorname{Spec }B} \rightarrow G(m,n)$ follows from a routine computation, which we omit.  Since ${\Phi_{i}}_{R}(\psi)$ is $A$-invariant, $\Phi_{i}$ factors through the inclusion $F_{A}(m,n) \rightarrow G(m,n)$.  If $K$ is infinite and $A$ is commutative, Corollary \ref{cor.before} implies that ${\Phi_{i}}_{R}(\psi)$ is generated by $A$-invariants.  Thus, $\Phi_{i}$ factors through the inclusion $H_{A}(m,n) \rightarrow G(m,n)$.
\end{proof}
We abuse notation by denoting both factors in the above lemma by $\Phi_{i}$.

\begin{lemma} \label{lemma.openo}
$\Phi_{i}:h_{\operatorname{Spec }B} \rightarrow F_{A}(m,n)$ is an open subfunctor.  Furthermore, if $K$ is infinite and $A$ is commutative, then $\Phi_{i}:h_{\operatorname{Spec }B} \rightarrow H_{A}(m,n)$ is an open subfunctor.
\end{lemma}

\begin{proof}
Suppose $\Psi:h_{\operatorname{Spec }R} \rightarrow F_{A}(m,n)$ is a natural transformation.  By Lemma \ref{lemma.welldef}, we must prove that, if $h_{\operatorname{Spec }B,\Psi}$ is the pullback in the diagram
$$
\begin{CD}
h_{\operatorname{Spec }B,\Psi} &  & \overset{\Gamma}{\longrightarrow} & & h_{\operatorname{Spec }R} \\
@VVV & & @VV{\Psi}V \\
h_{\operatorname{Spec }B} & & \underset{\Phi_{i}}{\longrightarrow} & & F_{A}(m,n)
\end{CD}
$$
then the induced map $\Gamma:h_{\operatorname{Spec }B,\Psi} \rightarrow h_{\operatorname{Spec }R}$ corresponds to the inclusion of an affine open subscheme of $\operatorname{Spec }R$.

By \cite[Exercise VI-6, p. 254]{eisen}, this is equivalent to showing that there exists some ideal $I$ of $R$ such that, for any $K$-algebra $T$,
\begin{equation} \label{eqn.open}
\operatorname{im}\Gamma_{T}=\{\delta \in h_{\operatorname{Spec }R}(T)|\delta(I)T=T\}.
\end{equation}
Let $f_{1},\ldots, f_{m} \in R^{lm}$ denote a basis for $\Psi(\operatorname{id}_{R})$, and suppose $f_{j}$ has $i$th coordinate $f_{ij}$.  Let $a$ denote the $m \times m$-matrix whose $p$th column is $(f_{(i-1)m+1,p},\ldots,f_{im,p})^{t}$, and let $I=\langle \operatorname{det }a \rangle$.  We prove that $I$ satisfies (\ref{eqn.open}).  That is, we prove that a homomorphism $\delta:R \rightarrow T$ has the property that
\begin{equation} \label{eqn.rightform}
\Psi(\delta)=\langle (\beta(x_{1}),\ldots,\beta(x_{(i-1)m}),1,0,\ldots,0,\beta(x_{(i-1)m+1}),\ldots,\beta(x_{lm-m})) \rangle
\end{equation}
for some $\beta:B \rightarrow T$ iff $\delta(I)T=T$.

Since $I$ is principle, $\delta:R \rightarrow T$ is such that $\delta(I)T=T$ iff $\delta(\operatorname{det }a)$ is a unit in $T$, which occurs iff the $m \times m$-matrix whose $p$th column is $(\delta(f_{(i-1)m+1,p}),\ldots,\delta(f_{im,p}))^{t}$ is invertible.  By naturality of $\Psi$, $\Psi(\delta)$ is the image of the composition
$$
T \otimes_{R} \Psi(\operatorname{id}_{R}) \rightarrow T \otimes_{R} R^{n} \rightarrow T^{n}
$$
whose left arrow is induced by inclusion $\Psi(\operatorname{id}_{R}) \rightarrow R^{n}$.  Thus, if $\delta(f_{j})$ denotes $(\delta(f_{1j}),\ldots,\delta(f_{lm,j}))^{t} \in T^{lm}$, then $\delta(f_{1}),\ldots,\delta(f_{m})$ is a basis for $\Psi(\delta)$.  This implies that the $m \times m$-matrix whose $p$th column is $(\delta(f_{(i-1)m+1,p}),\ldots,\delta(f_{im,p}))^{t}$ is invertible iff the projection of $\Psi(\delta)$ to the $(i-1)m+1$ through the $im$th factors is onto.  This occurs iff
$$
N=\langle(\beta(x_{1}),\ldots,\beta(x_{(i-1)m}),1,0,\ldots,0,\beta(x_{(i-1)m+1}),\ldots,\beta(x_{lm-m})) \rangle \subset \Psi(\delta)
$$
for some $\beta:B \rightarrow T$.  By Lemma \ref{lemma.rankm}, $N$ is a free rank $m$ summand of $R^{ml}$.  We claim $N=\Psi(\delta)$.  For, if $\mathfrak{m}$ is a maximal ideal of $T$, it follows from Nakayama's Lemma \cite[Corollary 4.8, p. 124]{eisenbud} that $N_{\mathfrak{m}}=\Psi(\delta)_{\mathfrak{m}}$.  Hence, $N=\Psi(\delta)$, and the first assertion follows.  To prove the second assertion, we note that the previous argument holds, mutatis mutandis, after replacing $F_{A}(m,n)$ by $H_{A}(m,n)$.
\end{proof}

\begin{cor} \label{lemma.open}
The open subfunctors $\Phi_{i}$ of $F_{A}(m,n)$ cover the $K$-rational points of $F_{A}(m,n)$.  That is,
$$
F_{A}(m,n)(K)=\underset{i}{\cup}\Phi_{i}(h_{\operatorname{Spec }B}(K)).
$$
Furthermore, if $K$ is infinite and $A$ is commutative, the open subfunctors $\Phi_{i}$ of $H_{A}(m,n)$ cover the $K$-rational points of $H_{A}(m,n)$.
\end{cor}

\begin{proof}
By Lemma \ref{lemma.openo}, the functors $\Phi_{i}:h_{\operatorname{Spec }B} \rightarrow F_{A}(m,n)$ are open.  If $M$ is a free rank $m$ summand of $K^{ml}$ which is $A$-invariant, then there exists an $i$ such that some element of $M$ has nonzero projection onto the $(i-1)m+1$ through the $im$th coordinates.  Hence, since $S$ is simple, $M$ contains the submodule
$$
\langle (b_{1},\ldots,b_{(i-1)m},1,0,\ldots,0,b_{(i-1)m+1},\ldots,b_{lm-m}) \rangle \subset K^{lm},
$$
where $b_{1},\ldots,b_{lm-m} \in K$.  Thus, by Lemma \ref{lemma.rankm},
$$
M=\langle (b_{1},\ldots,b_{(i-1)m},1,0,\ldots,0,b_{(i-1)m+1},\ldots,b_{lm-m}) \rangle,
$$
so that $M \in \Phi_{i}(h_{\operatorname{Spec }B}(K))$.  To prove the second assertion, we note that when $K$ is infinite and $A$ is commutative, Lemma \ref{lemma.openo} implies that $\Phi_{i}$ is an open subfunctor of $H_{A}(m,n)$.  Thus, the second assertion follows from the first assertion and Remark \ref{remark.newpoints}.
\end{proof}
We will use the following lemma to prove Theorem \ref{theorem.irr}.  The proof of the lemma is straightforward, so we omit it.

\begin{lemma} \label{lemma.irred}
Let $X$ be a topological space with open cover $\{A_{i}\}_{i \in I}$.  If $A_{i}$ is irreducible for all $i$ and $A_{i}\cap A_{j}$ is nonempty for all $i,j \in I$, then $X$ is irreducible.
\end{lemma}

\begin{thm} \label{theorem.irr}
Let $\mathbb{F}$ denote the open subscheme of $\mathbb{F}_{A}(m,n)$ obtained by glueing the open subschemes of $\mathbb{F}_{A}(m,n)$ defined by $\Phi_{i}$ for $1 \leq i \leq l$.  Then $\mathbb{F}$ is smooth, reduced, irreducible, of dimension $lm-m$ and has the same $K$-rational points as $\mathbb{F}_{A}(m,n)$.  If $K$ is infinite and $A$ is commutative, $\mathbb{H}_{A}(m,n)$ contains a smooth, reduced, irreducible, open subscheme of dimension $lm-m$ which has the same $K$-rational points as $\mathbb{H}_{A}(m,n)$.
\end{thm}

\begin{proof}
The fact that $\mathbb{A}_{K}^{lm-m}$ is smooth, reduced and has dimension $lm-m$ implies that $\mathbb{F}$ is smooth, reduced and has dimension $lm-m$.  The fact that $\mathbb{F}_{A}(m,n)$ and $\mathbb{F}$ have the same $K$-rational points follows from Corollary \ref{lemma.open}.

We denote the topological space of the open subscheme of $\mathbb{F}$ corresponding to $\Phi_{i}$ by $\mathbb{A}_{i}$.  For any $i,j$, the set $\mathbb{A}_{i} \cap \mathbb{A}_{j}$ is nonempty.  For example, if $i<j$, the intersection contains the point $(0,\ldots,0,1,0, \ldots, 0) \in \mathbb{A}_{i}$, where the nonzero entry occurs in the $(j-2)m+1$ position.  The fact that $\mathbb{F}$ is irreducible now follows from Lemma \ref{lemma.irred}.  The proof of the second assertion is similar, and we omit it.
\end{proof}

\section{Two-sided vector spaces}
In this section, we describe our notation and conventions regarding two-sided vector spaces, and we define the notion of rank of a two-sided vector space.  We end the section by reviewing facts about simple two-sided vector spaces which are employed in the sequel.

Let $V$ be a two-sided vector space.  That is, $V$ is a $k$-central $K-K$-bimodule which is finite-dimensional as a left $K$-module.  Right multiplication by $x\in K$ defines an endomorphism $\phi(x)$ of $_KV$, and the right action of $K$ on $V$ is via the $k$-algebra homomorphism $\phi:K\rightarrow \End(_KV)$.  This motivates the following definition.

\begin{defn} Let $\phi:K\rightarrow M_n(K)$ be a nonzero homomorphism.  We denote by $K^n_\phi$ the two-sided vector space of left dimension $n$, where the left action is the usual one and the right action is via $\phi$; that is,
\begin{equation}x\cdot(v_1,\dots, v_n)=(xv_1,\dots,xv_n),\ \ \ (v_1,\dots, v_n)\cdot x=(v_1,\dots,v_n)\phi(x).\end{equation}
We shall always write scalars as acting to the left of elements of $K^n_\phi$ and matrices acting to the right.
\end{defn}
If $V$ is a two-sided vector space, we let $\operatorname{dim }_{K}V$ denote the dimension of $V$ as a left $K$-module. If $\operatorname{dim }_{K}V=n$, then choosing a left basis for $V$ shows that $V\cong {K^n_\phi}$ for some homomorphism $\phi:K\rightarrow M_n(K)$.  Throughout the rest of this paper, $V$ will denote the two-sided vector space $K^n_\phi$, $W$ will denote a two-sided vector space and $S$ will denote a simple two-sided vector space.

We denote the category of two-sided vector spaces by $\Vect K$.  We shall denote by $K^B_i(K)$ the Quillen $K$-theory of $\Vect(K)$ (the superscript stands for ``bimodule").  The groups $K^{B}_{i}(B)$ were computed in \cite[Theorem 4.1]{pappacena}.

\begin{defn}
The rank of a two-sided vector space $W$, denoted $[W]$, is the class of $W$ in $K_{0}^{B}(K)$.
\end{defn}
Thus, the rank of $W$ is just the sums of the ranks of the simples (with multiplicity) appearing in the composition series of $W$.

We conclude this section with a description of the simple objects in $\Vect K$.  Let $\overline{K}$ denote an algebraic closure of $K$.  We write $\Emb(K)$ for the set of $k$-linear embeddings of $K$ into $\bar K$, and $G=G(K)$ for the group $\Aut(\bar K/K)$.

The group $G$ acts on $\Emb(K)$ by left composition. Given $\lambda\in \Emb(K)$, we denote the orbit of $\lambda$ under this action by $\lambda^G$, and we write $K(\lambda)$ for the composite field $K\vee\im(\lambda)$.

We denote the set of finite orbits of $\Emb(K)$ under the action of $G$ by $\Lambda(K)$.  The following is a consequence of the proof of \cite[Theorem 3.2]{pappacena}:

\begin{thm}  \label{theorem.main} If $K$ is perfect, there is a bijection from simple objects in $\Vect(K)$ to $\Lambda(K)$. Moreover, if $V$ is a simple two-sided vector space mapping to $\lambda^G\in \Lambda(K)$, and if $\lambda^G=\{\sigma_{1}\lambda, \ldots, \sigma_{m}\lambda \}$, then $\dim _KV =|\lambda^G|$ and there is a basis for the image of the composition
$$
K(\lambda) \otimes_{K} V \overset{=}{\rightarrow} K(\lambda) \otimes_{K} K^{n} \overset{\cong}{\rightarrow} K(\lambda)^{n}
$$
in which $\phi$ is a diagonal matrix with entries $\sigma_{1}\lambda, \ldots, \sigma_{m}\lambda$.
\end{thm}
We denote the simple two-sided vector space corresponding to $\lambda^{G}$ under the bijection in Theorem \ref{theorem.main} by $V(\lambda)$.

We will need the following Corollary to \cite[Lemma 3.13]{smith}:
\begin{lemma} \label{lemma.smith}
Let $F$ denote an extension field of $k$.  If $S$ and $S'$ are left finite-dimensional, non-isomorphic simple $F \otimes_{k}K$-modules, then $\operatorname{Ext}_{F \otimes_{k}K}^{1}(S,S')=0$.
\end{lemma}
Since a two-sided vector space is just a $K \otimes_{k} K$-module, Lemma \ref{lemma.smith} implies that $V \cong V_{1} \oplus \cdots \oplus V_{r}$, where $V_{i}$ is $S_{i}$-homogeneous for some simple $S_{i}$.

\section{Parameter spaces of two-sided subspaces of $V$}
The purpose of this section is to use $\mathbb{F}_{A}(m,n)$, $\mathbb{G}_{A}(m,n)$, and $\mathbb{H}_{A}(m,n)$ to construct and study parameter spaces of two-sided subspaces of $V$.

\subsection{The functors $F_{\phi}([W],V)$, $G_{\phi}([W],V)$, and $H_{\phi}([W],V)$}
\begin{defn}
If $V$ is $S$-homogeneous and $W$ is a two-sided vector space of rank $q[S]$, we let $F_{\phi}([W],V)(-): K-{\sf alg} \rightarrow {\sf Sets}$ denote the functor $F_{\operatorname{im }\phi}(qm,n)$.

If $W$ is not $S$-homogeneous, we let $F_{\phi}([W],V)(R)=\emptyset$.

Now suppose $V=V_{1} \oplus \cdots \oplus V_{r}$, where $V_{i}$ is $S_{i}$-homogeneous and $S_{i}$ is simple, $\phi_{i}(x)$ is the restriction of $\phi(x)$ to $V_{i}$, and $[W]=q_{1}[S_{1}]+\cdots+q_{r}[S_{r}]$.  We let $F_{\phi}([W],V)(-): K-{\sf alg} \rightarrow {\sf Sets}$ denote the functor
$$
F_{\phi_{1}}(q_{1}[S_{1}],V_{1}) \times \cdots \times F_{\phi_{r}}(q_{r}[S_{r}],V_{r}),
$$
where the product is taken over $h_{\operatorname{Spec }K}$.

If $W$ has a composition factor not in $\{S_{1},\ldots,S_{r}\}$, we let $F_{\phi}([W],V)(R)=\emptyset$.
\end{defn}
We call elements of $F_{\phi}([W],V)(R)$ {\it free rank $[W]$ $\phi$-invariant families over $\operatorname{Spec }R$}, or {\it free $\phi$-invariant families} when $W$ and $R$ are understood.

\begin{defn}
If $V$ is $S$-homogeneous and $W$ is a two-sided vector space of rank $q[S]$, we let $G_{\phi}([W],V)(-): K-{\sf alg} \rightarrow {\sf Sets}$ denote the functor $G_{\operatorname{im }\phi}(qm,n)$.

If $W$ is not $S$-homogeneous, we let $G_{\phi}([W],V)(R)=\emptyset$.

Now suppose $V=V_{1} \oplus \cdots \oplus V_{r}$, where $V_{i}$ is $S_{i}$-homogeneous and $S_{i}$ is simple, $\phi_{i}(x)$ is the restriction of $\phi(x)$ to $V_{i}$, and $[W]=q_{1}[S_{1}]+\cdots+q_{r}[S_{r}]$.  We let $G_{\phi}([W],V)(-): K-{\sf alg} \rightarrow {\sf Sets}$ denote the functor
$$
G_{\phi_{1}}(q_{1}[S_{1}],V_{1}) \times \cdots \times G_{\phi_{r}}(q_{r}[S_{r}],V_{r}),
$$
where the product is taken over $h_{\operatorname{Spec }K}$.

If $W$ has a composition factor not in $\{S_{1},\ldots,S_{r}\}$, we let $G_{\phi}([W],V)(R)=\emptyset$.
\end{defn}
We call elements of $G_{\phi}([W],V)(R)$ {\it free rank $[W]$ families generated by $\phi$-invariants over $\operatorname{Spec }R$}, or {\it free families generated by $\phi$-invariants} when $W$ and $R$ are understood.

\begin{defn}
If $V$ is $S$-homogeneous and $W$ is a two-sided vector space of rank $q[S]$, we let $H_{\phi}([W],V)(-): K-{\sf alg} \rightarrow {\sf Sets}$ denote the functor $H_{\operatorname{im }\phi}(qm,n)$.

If $W$ is not $S$-homogeneous, we let $H_{\phi}([W],V)(R)=\emptyset$.

Now suppose $V=V_{1} \oplus \cdots \oplus V_{r}$, where $V_{i}$ is $S_{i}$-homogeneous and $S_{i}$ is simple, $\phi_{i}(x)$ is the restriction of $\phi(x)$ to $V_{i}$, and $[W]=q_{1}[S_{1}]+\cdots+q_{r}[S_{r}]$.  We let $H_{\phi}([W],V)(-): K-{\sf alg} \rightarrow {\sf Sets}$ denote the functor
$$
H_{\phi_{1}}(q_{1}[S_{1}],V_{1}) \times \cdots \times H_{\phi_{r}}(q_{r}[S_{r}],V_{r}),
$$
where the product is taken over $h_{\operatorname{Spec }K}$.

If $W$ has a composition factor not in $\{S_{1},\ldots,S_{r}\}$, we let $H_{\phi}([W],V)(R)=\emptyset$.
\end{defn}
We call elements of $H_{\phi}([W],V)(R)$ {\it free rank $[W]$ $\phi$-invariant families generated by $\phi$-invariants over $\operatorname{Spec }R$}, or {\it free $\phi$-invariant families generated by $\phi$-invariants} when $W$ and $R$ are understood.

\begin{lemma}
The $K$-rational points of $F_{\phi}([W],V)$, $G_{\phi}([W],V)$, and $H_{\phi}([W],V)$ are equal to the set of two-sided rank $[W]$ subspaces of $V$.
\end{lemma}

\begin{proof}
We first show that the three functors above have the same $K$-rational points.  From the definitions of $F_{\phi}([W],V)$, $G_{\phi}([W],V)$, and $H_{\phi}([W],V)$, it suffices to prove the result when $V$ is homogeneous.  Thus, it suffices to prove
$$
F_{\operatorname{im }\phi}(m,n)(K)=G_{\operatorname{im }\phi}(m,n)(K)=H_{\operatorname{im }\phi}(m,n)(K)
$$
when $K^{n}$ is homogeneous as a $K\otimes_{k}\operatorname{im }\phi$-module.  Since $\phi:K \rightarrow M_{n}(K)$ is a ring homomorphism, $\operatorname{im }\phi -\{0\} \subset GL_{n}(K)$.  Thus, the assertion follows from Remark \ref{remark.points} and Remark \ref{remark.newpoints}.

To complete the proof of the lemma, it suffices to prove that
$$
F_{\phi}([W],V)(K)=\{\mbox{two-sided rank $[W]$ subspaces of $V$}\}.
$$
If $V$ is homogeneous, this follows immediately from the definition of $F_{\operatorname{im }\phi}(m,n)$.  If $V$ is not homogeneous, the result follows from the fact that $V$, and any two-sided subspace of $V$, has a direct sum decomposition into its homogeneous components.
\end{proof}

We now find conditions under which $F_{\phi}([S],V) \neq G_{\phi}([S],V)$ and $F_{\phi}([S],V) \neq H_{\phi}([S],V)$.

\begin{lemma} \label{lemma.arturo}
Suppose $\lambda_{1},\ldots,\lambda_{m} \in \operatorname{Emb}(K)$ are distinct and $|k|>m$.  If
$$
\{i_{1},\ldots,i_{m} \} \subset \{1, \ldots, m\}
$$
is a multiset with repetitions, then there exists an $a \in K$ such that
\begin{equation} \label{eqn.newstar}
\prod_{j=1}^{m}\lambda_{j} (a) \neq \prod_{i=1}^{m}\lambda_{i_{j}}(a).
\end{equation}
\end{lemma}

\begin{proof}
First we claim that there exists an element $b \in K$ such that there is an inequality of multisets
$$
\{\lambda_{1}(b),\ldots,\lambda_{m}(b)\} \neq \{\lambda_{i_{1}}(b),\ldots,\lambda_{i_{m}}(b)\}.
$$
If not, we would have
$$
\sum_{j=1}^{m} \lambda_{j} = \sum_{j=1}^{m} \lambda_{i_{j}},
$$
which is a nontrivial dependency relation among $\{\lambda_{1}, \ldots, \lambda_{m}\}$ as $k$-linear functions from $K$ to $\overline{K}$.  This contradicts the linear independence of characters from $K$ to $\overline{K}$, which establishes our claim.

With $b$ as above, we have
$$
\prod_{j=1}^{m}(x-\lambda_{j}(b)) \neq \prod_{j=1}^{m}(x-\lambda_{i_{j}}(b))
$$
in the ring $\overline{K}[x]$.  Thus
$$
f(x)=\prod_{j=1}^{m}(x-\lambda_{j}(b)) - \prod_{j=1}^{m}(x-\lambda_{i_{j}}(b))
$$
has at most $m$ roots.  Since $|k|>m$, we may choose $c \in k$ such that $f(c) \neq 0$.  Since $\lambda_{i}$ is $k$-linear, we have
$$
\prod_{j=1}^{m}(\lambda_{j}(c-b)) - \prod_{j=1}^{m}(\lambda_{i_{j}}(c-b)) = \prod_{j=1}^{m}(c-\lambda_{j}(b)) - \prod_{j=1}^{m}(c-\lambda_{i_{j}}(b)) = f(c) \neq 0.
$$
Thus (\ref{eqn.newstar}) holds with $a=c-b$.
\end{proof}

\begin{cor} \label{cor.noninvar}
Suppose $m \in \mathbb{N}$ is such that $|k|>m>1$, $K$ is perfect and $\lambda \in \operatorname{Emb }(K)$ is such that $|\lambda^{G}|=m$ (see Section 5 for notation).  If $V(\lambda)^{\oplus 2} \subset V$, then there exist free rank $[V(\lambda)]$ $\phi$-invariant families over $\operatorname{Spec }K(\lambda)$ which are not generated by $\operatorname{im }\phi$-invariants.
\end{cor}

\begin{proof}
Suppose $\lambda^{G}=\{\sigma_{1}\lambda, \ldots, \sigma_{m} \lambda\}$.  Let $M$ be a free rank $[V(\lambda)]$ $\phi$-invariant family over $\operatorname{Spec }K(\lambda)$ which is generated by $\phi$-eigenvectors $v_{1},\ldots,v_{m}$, with eigenvalues
$$
\sigma_{i_{1}}\lambda,\ldots,\sigma_{i_{m}}\lambda
$$
respectively, such that the multiset $\{i_{1},\ldots,i_{m}\}$ has repetitions.  Then the eigenvalue of any generator of $\bigwedge^{m}M$ equals $\prod_{j}\sigma_{i_{j}}\lambda$.  On the other hand, if $M$ were generated by $\operatorname{im }\phi$-invariants, any generator of $\bigwedge^{m}M$ would have eigenvalue $\prod_{j}\sigma_{j}\lambda$.  Thus, by Lemma \ref{lemma.arturo}, $M$ is a free rank $[V(\lambda)]$ $\phi$-invariant family which is not generated by $\operatorname{im }\phi$-invariants.
\end{proof}

\begin{example} \label{example.notthesame}
Suppose $\rho = {}^{3}\sqrt{2}$, $\zeta$ is a primitive $3$rd root of unity, $k=\mathbb{Q}$ and $K=\mathbb{Q}(\rho)$.  For $i=0,1$, let
$$
\lambda_{i}(\sum_{l=0}^{2}a_{l}\rho^{l})=a_{i}\rho^{i}-a_{2}\rho^{2}
$$
and let $\lambda(x)=\lambda_{0}(x)+\lambda_{1}(x)\zeta$.  Then $V(\lambda)$ is a two-dimensional simple two-sided vector space \cite[Example 3.9]{pappacena}, and thus, by Corollary \ref{cor.noninvar}, $V=V(\lambda)^{\oplus 2}$ contains free rank $V(\lambda)$ $\phi$-invariant families over $\operatorname{Spec }K(\lambda)$ which are not generated by $\operatorname{im }\phi$-invariants.  In other words,
$$
F_{\phi}([V(\lambda)],V(\lambda)^{\oplus 2})(K(\lambda)) \neq H_{\phi}([V(\lambda)],V(\lambda)^{\oplus 2})(K(\lambda)).
$$
It follows immediately that
$$
F_{\phi}([V(\lambda)],V(\lambda)^{\oplus 2})(K(\lambda)) \neq G_{\phi}([V(\lambda)],V(\lambda)^{\oplus 2})(K(\lambda)).
$$
\end{example}

\begin{remark}
It follows from the previous example and the definitions of $F_{\phi}([W],V)$, $G_{\phi}([W],V)$, and $H_{\phi}([W],V)$ that there exist $A$, $m$ and $n$ such that $F_{A}(m,n) \neq G_{A}(m,n)$ and $F_{A}(m,n) \neq H_{A}(m,n)$.
\end{remark}

\subsection{$F$-rational points of $G_{\phi}([S],V)$ and $H_{\phi}([S],V)$}
Let $F$ be an extension field of $K$.  In this subsection, we show that every element of $G_{\phi}([S],V)(F)$ and of $H_{\phi}([S],V)(F)$ is isomorphic to $F \otimes_{K} S$ as $F\otimes_{k}K$-modules.  Throughout this subsection, we assume, without loss of generality, that $V$ is $S$-homogeneous.  Since, by Lemma \ref{lemma.phiphiinvar}, $G_{\phi}([S],V)(F)=H_{\phi}([S],V)(F)$, it suffices to prove the result for $G_{\phi}([S],V)(F)$.  We assume throughout this subsection that $\operatorname{dim }_{K}S=m$.

Since the sum of all simple submodules of $V$ is a direct summand of $V$ as a left $K$-module, $V$ has a left $K$-module decomposition
\begin{equation} \label{eqn.decomp}
V=L \oplus N,
\end{equation}
where $N$ contains no simple two-sided subspaces of $V$ and
$L=S^{\oplus l}$ is a direct sum of simple two-sided subspaces of
$V$.
\begin{lemma} \label{lemma.proj}
Every free rank $[S]$ family generated by $\phi$-invariants over $\operatorname{Spec }F$ is contained in $F \otimes_{K} L$.
\end{lemma}

\begin{proof}
Assume $N \neq 0$ and suppose $M$ is a free rank $[S]$ family generated by $\phi$-invariants over $\operatorname{Spec }F$, with basis $\{v_{i}+w_{i}\}_{i=1}^{m}$, where $v_{i}$ is an element of the image of the composition
\begin{equation} \label{eqn.image}
F \otimes_{K}N \rightarrow F \otimes_{K}V \overset{\cong}{\rightarrow} F^{n}
\end{equation}
induced by the inclusion $N \subset V$, and $w_{i}$ is an element of the image of the composition
\begin{equation} \label{eqn.imagea}
F \otimes_{K} L \rightarrow F \otimes_{K} V \overset{\cong}{\rightarrow} F^{n}
\end{equation}
induced by the inclusion $L \subset V$.  Since $M$ is generated by $\operatorname{im }\phi$-invariants, $\wedge^{m} M$ is contained in the image of the composition
$$
F \otimes_{K}\textstyle\bigwedge^{m}L \rightarrow F \otimes_{K}\textstyle\bigwedge^{m}V \rightarrow \textstyle\bigwedge^{m}F^{n}
$$
induced by the inclusion $\bigwedge^{m}L \rightarrow \bigwedge^{m}V$.  On the other hand,
\begin{equation} \label{eqn.expand}
(v_{1}+w_{1})\wedge \cdots \wedge (v_{m}+w_{m})=w_{1}\wedge \cdots \wedge w_{m}+\mbox{wedges with at least one $v_{i}$.}
\end{equation}
Since $\bigwedge^{m}M \neq 0$, this implies that $w_{1}\wedge \cdots \wedge w_{m} \neq 0$, and that the sum of the other terms on the right-hand side of (\ref{eqn.expand}) equals $0$.  We claim each of the terms $v_{1}\wedge w_{2} \wedge \cdots \wedge w_{m}, \ldots, w_{1}\wedge \cdots \wedge w_{m-1}\wedge v_{m}$ equals zero, which would prove the assertion.  To this end, let $f_{1},\ldots, f_{p}$ denote a basis for the image of (\ref{eqn.image}), and suppose $f_{p+1},\ldots,f_{q}, w_{1}, \ldots, w_{m}$ is a basis for the image of (\ref{eqn.imagea}).  Let
$$
B_{1}=\{w_{1}\wedge \cdots \wedge w_{m}\} \cup \bigg{(}\underset{1 \leq j_{1} < \cdots < j_{m} \leq q}{\bigcup}\{f_{j_{1}}\wedge \cdots \wedge f_{j_{m}}\}\bigg{)},
$$
$$
B_{2}=\overset{m-1}{\underset{r=2}{\bigcup}}\bigg{(} \underset{1 \leq i_{r+1} < \cdots < i_{m} \leq m}{\underset{1 \leq i_{1} < \cdots < i_{r} \leq q}{\bigcup}} \{f_{i_{1}}\wedge \cdots \wedge f_{i_{r}} \wedge w_{i_{r+1}} \wedge \cdots \wedge w_{i_{m}}\} \bigg{)},
$$
$$
B_{3}= \bigg{(}\overset{q}{\underset{i=1}{\bigcup}}\{w_{1}\wedge f_{i} \wedge w_{3} \wedge \cdots \wedge w_{m}\}\bigg{)} \cup \cdots \cup \bigg{(}\overset{q}{\underset{i=1}{\bigcup}}\{w_{1}\wedge \cdots \wedge w_{m-1}\wedge f_{i}\}\bigg{)},
$$
and
$$
B_{4}=\overset{q}{\underset{i=1}{\bigcup}}\{f_{i}\wedge w_{2} \wedge w_{3} \wedge \cdots \wedge w_{m}\}.
$$
The sets $B_{1},B_{2},B_{3},B_{4}$ form a partition of a basis for $\bigwedge^{m}F^{n}$.  Since the right-hand side of (\ref{eqn.expand}) equals $w_{1}\wedge \cdots \wedge w_{m}$, $v_{1} \wedge w_{2} \wedge \cdots \wedge w_{m}$ is a linear combination of elements in $B_{1} \cup B_{2} \cup B_{3}$.  On the other hand, $v_{1} \wedge w_{2} \wedge \cdots \wedge w_{m}$ is a linear combination of elements in $B_{4}$.  We conclude that $v_{1} \wedge w_{2} \wedge \cdots \wedge w_{m}=0$.  A similar argument implies that $w_{1}\wedge \cdots \wedge w_{i-1} \wedge v_{i} \wedge w_{i+1} \wedge \cdots \wedge w_{m}=0$ for $1 \leq i \leq m$, and the result follows.
\end{proof}

\begin{thm} \label{theorem.proj}
Suppose $|k|>m$, $K$ is perfect, and $K \subset F$ is an extension of fields.  If $M$ is a free rank $[S]$ family generated by $\phi$-invariants over $\operatorname{Spec }F$, then $M \cong F \otimes_{K}S$ as $F \otimes_{k}K$-modules.
\end{thm}

\begin{proof}
By Lemma \ref{lemma.proj}, we may assume $M$ is contained in the image of the composition
$$
F \otimes_{K} L \rightarrow F \otimes_{K}V \overset{\cong}{\rightarrow} F^{n}
$$
induced by the inclusion $L \subset V$.  Thus, we may assume $V$ is semisimple.

Let $\overline{F}$ denote an algebraic closure of $F$ containing $\overline{K}$, let $\overline{M}$ denote the image of the composition
$$
\overline{F} \otimes_{F} M \rightarrow \overline{F} \otimes_{F} F^{n} \overset{\cong}{\rightarrow} \overline{F}^{n}
$$
whose left arrow is induced by inclusion, and let $\overline{M}$ have generators $w_{1},\ldots, w_{m}$ as an $\overline{F}$-module.  By Lemma \ref{lemma.functor}, $\overline{M}$ is generated by $\operatorname{im }\phi$-invariants.  Thus, Theorem \ref{theorem.main} implies that the $\phi$-eigenvalues of $w_{1}\wedge \cdots \wedge w_{m}$ must equal $\sigma_{1}\lambda (x) \cdots \sigma_{m}\lambda (x)$ for all $x \in K$, where $\lambda$ is a $k$-linear embedding of $K$ into $\overline{K}$, $\sigma_{1},\ldots,\sigma_{m}$ are automorphisms of $\overline{K}$ over $K$, and $\{\sigma_{1} \lambda, \ldots, \sigma_{m} \lambda\}$ are distinct.  By Lemma \ref{lemma.phiphiinvar}, $\overline{M}$ is $\phi$-invariant.  Thus, $\overline{M}$ has a $\phi$-eigenvector, $v_{1}$.  Since $v_{1}$ is also an eigenvector in $\overline{F}^{n} \supset \overline{K}^{n}$, it must have eigenvalue $\sigma_{i_{1}}\lambda$.  For $1 < j \leq m$, let $v_{j}\in \overline{M}$ be such that $v_{j}+\langle v_{1},\ldots,v_{j-1} \rangle$ is a $\phi$-eigenvector for $\overline{M}/\langle v_{1}, v_{2}, \ldots,v_{j-1} \rangle$, where $\langle v_{1},\ldots,v_{j-1} \rangle$ denotes the $\overline{F}\otimes_{k}K$-module generated by $v_{1},\ldots,v_{j-1}$.  Then $v_{j}+\langle v_{1},\ldots,v_{j-1} \rangle$ has eigenvalue $\sigma_{i_{j}}\lambda$, and, thus, in the basis $\{v_{1},\ldots,v_{m}\}$, $\phi(x)|_{\overline{M}}$ is upper-triangular with diagonal entries $\sigma_{i_{1}}\lambda(x),\ldots,\sigma_{i_{m}}\lambda(x)$.  Therefore,
$$
\operatorname{det} \phi(x)|_{\overline{M}} = \sigma_{i_{1}}\lambda(x) \cdots \sigma_{i_{m}} \lambda(x).
$$
By Lemma \ref{lemma.arturo}, we must have $\{i_{1},\ldots, i_{m}\}=\{1, \ldots, m\}$.  Since, by Lemma \ref{lemma.smith}, extensions of distinct simple left finite-dimensional $\overline{F} \otimes_{k}K$-modules are split, there exists a basis for $\overline{M}$ such that $\phi(x)|_{\overline{M}}$ is diagonal with entries $\sigma_{1}\lambda(x),\ldots, \sigma_{m}\lambda(x)$.  It follows that $\overline{F} \otimes_{F}M \cong \overline{F} \otimes_{F} (F \otimes_{K} S)$ as $\overline{F} \otimes_{k} K$-modules.  Thus, by an argument similar to that given in the proof of \cite[Lemma 2.4]{pappacena}, we conclude that $M \cong F\otimes_{K}S$ as $F \otimes_{k} K$-modules.
\end{proof}

\subsection{The geometry of $\mathbb{F}_{\phi}([W],V)$, $\mathbb{G}_{\phi}([W],V)$, and $\mathbb{H}_{\phi}([W],V)$}
For the readers convenience, we collect here some consequences of our study of the geometry of $F_{A}(m,n)$, $G_{A}(m,n)$ and $H_{A}(m,n)$ in the case that $A = \operatorname{im }\phi$, where $\phi:K \rightarrow M_{n}(K)$ is a $k$-central ring homomorphism.  We assume throughout the remainder of this section that $[V]=l_{1}[S_{1}]+\cdots+l_{r}[S_{r}]$, where $S_{1},\ldots,S_{r}$ are non-isomorphic simple modules with $\operatorname{dim }S_{i}=m_{i}$, and that $\phi_{i}(x)$ is the restriction of $\phi(x)$ to the $S_{i}$-homogeneous summand of $V$.  Finally, we assume all products of schemes are over $\operatorname{Spec }K$.

\begin{thm}
The functors $F_{\phi}(q_{1}[S_{1}]+\cdots+q_{r}[S_{r}],V)$, $G_{\phi}(q_{1}[S_{1}]+\cdots+q_{r}[S_{r}],V)$, and $H_{\phi}(q_{1}[S_{1}]+\cdots+q_{r}[S_{r}],V)$ are represented by
$$
\overset{r}{\underset{i=1}{\prod}}\mathbb{F}_{\operatorname{im }\phi_{i}}(m_{i}q_{i},m_{i}l_{i}), \overset{r}{\underset{i=1}{\prod}}\mathbb{G}_{\operatorname{im }\phi_{i}}(m_{i}q_{i},m_{i}l_{i}), \mbox{ and }\overset{r}{\underset{i=1}{\prod}}\mathbb{H}_{\operatorname{im }\phi_{i}}(m_{i}q_{i},m_{i}l_{i})
$$
respectively.
\end{thm}

\begin{proof}
Since $F_{A}(m,n)$, $G_{A}(m,n)$, and $H_{A}(m,n)$ are representable by $\mathbb{F}_{A}(m,n)$, ${\mathbb{G}}_{A}(m,n)$, and ${\mathbb{H}}_{A}(m,n)$, the result follows from \cite[p. 260]{eisen}.
\end{proof}
We denote the schemes representing $F_{\phi}([W],V)$, $G_{\phi}([W],V)$, and $H_{\phi}([W],V)$ by $\mathbb{F}_{\phi}([W],V)$, $\mathbb{G}_{\phi}([W],V)$, and $\mathbb{H}_{\phi}([W],V)$, respectively.

\begin{cor} \label{cor.galois}
If $K/k$ is finite and Galois then $F_{\phi}([W],V)=G_{\phi}([W],V)=H_{\phi}([W],V)$ and ${\mathbb{F}}_{\phi}(q_{1}[S_{1}]+\cdots+q_{r}[S_{r}],V)$ equals
$$
\overset{r}{\underset{i=1}{\prod}}\mathbb{G}(q_{i},l_{i}).
$$
\end{cor}

\begin{proof}
We prove that ${\mathbb{F}}_{\phi}(q_{1}[S_{1}]+\cdots+q_{r}[S_{r}],V)=\overset{r}{\underset{i=1}{\prod}}\mathbb{G}(q_{i},l_{i})$.  The other assertions follow similarly.  By the previous result, it suffices to prove that $\mathbb{F}_{\operatorname{im }\phi_{i}}(m_{i}q_{i},m_{i}l_{i})=\mathbb{G}(q_{i},l_{i})$.  The hypothesis on $K/k$ is equivalent to $K$ being a finite, separable extension of $k$ such that $\operatorname{Aut} K=\operatorname{Emb} K$.  Thus, $m_{i}=1$ \cite[Theorem 3.2]{pappacena} and $S_{i} \cong K_{\sigma_{i}}$ for some $k$-linear automorphism $\sigma_{i}$ of $K$ (note that, in this case, we do not require that $K$ be perfect to apply \cite[Theorem 3.2]{pappacena}).  Since $K\otimes_{k}K$ is semisimple, $\phi_{i}$ is a diagonal matrix with each diagonal entry equal to $\sigma_{i}$, and the assertion follows.
\end{proof}

\begin{remark}
The previous result also follows from the second part of \cite[Theorem 1, p. 321]{pop}.
\end{remark}

Now assume that $V$ is semisimple.  Before we state our next result, we need to introduce some notation.  For $1 \leq i \leq r$, let $B_{i}=K[x_{i,1},\ldots, x_{i,l_{i}m_{i}-m_{i}}]$, let $B=K[\{x_{i,1},\ldots,x_{i,l_{i}m_{i}-m_{i}}\}_{i=1}^{r}]$ and, for each $r$-tuple $J=(j_{1},\ldots,j_{r})$ such that $1 \leq j_{i} \leq l_{i}$ define an inclusion of functors
$$
\Phi_{J}:h_{\operatorname{Spec }B_{1}} \times \cdots \times h_{\operatorname{Spec }B_{r}} \rightarrow F_{\phi_{1}}([S_{1}],V_{1}) \times \cdots \times F_{\phi_{r}}([S_{r}],V_{r})
$$
by $\Phi_{J}=\Phi_{j_{1}} \times \cdots \times \Phi_{j_{r}}$, where $\Phi_{i}$ is defined by (\ref{eqn.bigphi}), and where all products are over $h_{\operatorname{Spec }K}$.  We abuse notation by letting $\Phi_{J}$ denote the induced natural transformation
$$
h_{\operatorname{Spec }B} \overset{\cong}{\rightarrow} h_{\operatorname{Spec }B_{1}} \times \cdots \times h_{\operatorname{Spec }B_{r}} \overset{\Phi_{J}}{\rightarrow} F_{\phi}([S_{1}]+\cdots+[S_{r}],V).
$$
In a similar fashion, we can define an inclusion of functors
$$
\Phi_{J}:h_{\operatorname{Spec }B_{1}} \times \cdots \times h_{\operatorname{Spec }B_{r}} \rightarrow H_{\phi_{1}}([S_{1}],V_{1}) \times \cdots \times H_{\phi_{r}}([S_{r}],V_{r}),
$$
where we have abused notation as in Section 4.

The following result is an immediate consequence of Lemma \ref{lemma.openo} and Corollary \ref{lemma.open}.
\begin{thm}
For all $r$-tuples $J=(j_{1},\ldots,j_{r})$ such that $1 \leq j_{i} \leq l_{i}$, $\Phi_{J}: h_{\operatorname{Spec }B} \rightarrow F_{\phi}([S_{1}]+\cdots+[S_{r}],V)$ is an open subfunctor, and the open subfunctors $\Phi_{J}$ cover the $K$-rational points of $F_{\phi}([S_{1}]+\cdots+[S_{r}],V)$.  Furthermore, if $K$ is infinite, the same result holds for $H_{\phi}([S_{1}]+\cdots+[S_{r}],V)$ in place of $F_{\phi}([S_{1}]+\cdots+[S_{r}],V)$.
\end{thm}

The following follows from the above result and from an argument similar to that used to prove Theorem \ref{theorem.irr}.
\begin{cor} \label{cor.geometry1}
${\mathbb{F}}_{\phi}([S_{1}]+\cdots+[S_{r}],V)$ and ${\mathbb{H}}_{\phi}([S_{1}]+\cdots+[S_{r}],V)$ contain smooth, reduced, irreducible open subschemes of dimension $\sum_{i=1}^{r}l_{i}m_{i}-m_{i}$ which cover their $K$-rational points.
\end{cor}

The following example illustrates the fact that the open subfunctors $\Phi_{J}$ do not always form an open cover of $F_{\phi}([S],S^{\oplus l})$ or $H_{\phi}([S],S^{\oplus l})$.
\begin{example}
Suppose $\rho = {}^{3}\sqrt{2}$, $\zeta$ is a primitive $3$rd root of unity, $k=\mathbb{Q}$ and $K=\mathbb{Q}(\rho)$.  For $i=0,1$, let
$$
\lambda_{i}(\sum_{l=0}^{2}a_{l}\rho^{l})=a_{i}\rho^{i}-a_{2}\rho^{2}
$$
and let $\lambda(x)=\lambda_{0}(x)+\lambda_{1}(x)\zeta$.  Let $V(\lambda)$ denote the corresponding two-dimensional simple $K \otimes_{k}K$-module, so that the right action of $K$ on $V(\lambda)$ is given by $\phi(x)=\begin{pmatrix} \lambda_{0}(x) & -\lambda_{1}(x) \\ \lambda_{1}(x) & -\lambda_{1}(x)+\lambda_{0}(x) \end{pmatrix}$ \cite[Example 3.9]{pappacena}.  Let $V=V(\lambda)^{\oplus 2}$, and let $\{e_{i}\}_{i=1}^{4}$ denote the standard unit vectors of $K(\zeta)^{4}$.  Then
$$
M=\operatorname{Span}_{K(\zeta)}\{e_{1}+\zeta e_{2}, e_{3}+\zeta^{2} e_{4}\} \subset K(\zeta)^{4} \cong K(\zeta) \otimes_{K} V
$$
is a free $\phi$-invariant rank $[V(\lambda)]$ family over $\operatorname{Spec }K(\zeta)$ whose projections onto the first and second coordinates of $K(\zeta)^{4}$, and onto the third and fourth coordinates of $K(\zeta)^{4}$, are not onto.  In particular, $M$ is not an element of $\overset{2}{\underset{J=1}{\cup}}\Phi_{J}(h_{\operatorname{Spec }A}(K(\zeta)))$, and hence by \cite[Exercise VI-II, p. 256]{eisen}, the open subfunctors $\Phi_{J}$ do not cover $F_{\phi}([V(\lambda)],V(\lambda)^{\oplus 2})$.

We claim that $M$ is generated by $\operatorname{im }\phi$-invariants, which would establish that the open subfunctors $\Phi_{J}$ do not cover $H_{\phi}([V(\lambda)],V(\lambda)^{\oplus 2})$.  To prove the claim, we first note that
\begin{equation} \label{equation.zeta}
(e_{1}+\zeta e_{2})\wedge (e_{3}+\zeta^{2} e_{4})=e_{1}\wedge e_{3}+\zeta^{2} e_{1}\wedge e_{4} + \zeta e_{2}\wedge e_{3}+e_{2} \wedge e_{4}.
\end{equation}
On the other hand, if we let $w_{1}=e_{1}$, $w_{2}=e_{3}$, $a_{1}=1$ and $a_{2}=\rho$, then
$$
w_{1}\phi(a_{1})\wedge w_{2}\phi(a_{2})+w_{2}\phi(a_{1})\wedge w_{1}\phi(a_{2}) = -\rho(e_{1}\wedge e_{4}-e_{2} \wedge e_{3})
$$
is an element of ${\bigwedge}_{\operatorname{im }\phi}^{2}$ by Proposition \ref{proposition.phiinvar}.  Similarly, if we let $w_{1}=e_{1}$, $w_{2}=e_{4}$, $a_{1}=1$ and $a_{2}=\rho$, then
$$
w_{1}\phi(a_{1})\wedge w_{2}\phi(a_{2})+w_{2}\phi(a_{1})\wedge w_{1}\phi(a_{2}) = \rho(e_{1}\wedge e_{3} - e_{1}\wedge e_{4}-e_{4}\wedge e_{2})
$$
is an element of ${\bigwedge}_{\operatorname{im }\phi}^{2}$ by Proposition \ref{proposition.phiinvar}.  It follows that (\ref{equation.zeta}) is an element of the image of $K(\zeta) \otimes_{K} {\bigwedge}_{\operatorname{im }\phi}^{2} \rightarrow K(\zeta) \otimes_{K} \bigwedge^{2}K^{4} \overset{\cong}{\rightarrow} \bigwedge^{2}K(\zeta)^{4}$, and hence that $M$ is generated by $\operatorname{im }\phi$-invariants.
\end{example}

\end{document}